\documentclass[11pt]{article}

\usepackage{verbatim,latexsym,amsfonts,amsmath,amssymb,graphicx,fancyhdr,hyperref,asymptote}

\newcommand{\EE}{\mbox{\bf E}\,}
\newcommand{\PP}{\mathbb{P}}

\newcommand{\R}{\mathbb{R}}
\newcommand{\C}{\mathbb{C}}

\newcommand{\HH}{\mathbb{H}}
\newcommand{\N}{\mathbb{N}}
\newcommand{\D}{\mathbb{D}}
\newcommand{\Z}{\mathbb{Z}}
\newcommand{\St}{\mathbb{S}}
\newcommand{\SA}{\mbox{\bf S}}
\newcommand{\HA}{\mbox{\bf H}}
\newcommand{\TT}{\mathbb{T}}
\newcommand{\A}{\mathbb{A}}

\newcommand{\pa}{\partial}

\newcommand{\F}{{\cal F}}

\newcommand{\no}{\noindent}
\newcommand{\rA}{\mbox{\bf r}}
\newcommand{\RA}{\mbox{\bf R}}

\def\til{\widetilde}
\def\ha{\widehat}
\def\sem{\setminus}
\def\lin{\overline}
\def\vphi{\varphi}

\def\luto{\stackrel{\rm l.u.}{\longrightarrow}}
\def\r{\mathring}

\DeclareMathOperator{\ccap}{cap}

 \DeclareMathOperator{\diam}{diam}
\DeclareMathOperator{\dist}{dist} 
\DeclareMathOperator{\hcap}{hcap} \DeclareMathOperator{\id}{id}
\DeclareMathOperator{\Imm}{Im } \DeclareMathOperator{\Ree}{Re }

\DeclareMathOperator{\modd}{mod} 
\DeclareMathOperator{\PV}{P.V.}

 \DeclareMathOperator{\cc}{c}
\DeclareMathOperator{\disj}{disj}

\newtheorem{Lemma}{Lemma}[section]
\newtheorem{Theorem}{Theorem}[section]
\newtheorem{Definition}{Definition}[section]

\newtheorem{Proposition}{Proposition}[section]
\newtheorem{Example}{Example}
\numberwithin{equation}{section}
\newcommand{\BGE}{\begin{equation}}
\newcommand{\BGEN}{\begin{equation*}}
\newcommand{\EDE}{\end{equation}}
\newcommand{\EDEN}{\end{equation*}}
\def\dto{\stackrel{\rm Cara}{\longrightarrow}}
\def\conf{\stackrel{\rm Conf}{\twoheadrightarrow}}

\begin{document}

\title{\bf Restriction Properties of Annulus SLE}
\date{\today}
\author{Dapeng Zhan }
\maketitle

\begin{abstract}
For $\kappa\in(0,4]$, a family of annulus SLE$(\kappa;\Lambda)$ processes were introduced in \cite{whole} to prove the reversibility of whole-plane SLE$(\kappa)$. In this paper we prove that those annulus SLE$(\kappa;\Lambda)$ processes satisfy a restriction property, which is similar to that for chordal SLE$(\kappa)$. Using this property, we construct $n\ge 2$ curves crossing an annulus such that, when any $n-1$ curves are given, the last curve is a chordal SLE$(\kappa)$ trace.
\end{abstract}


\section{Introduction}
Oded Schramm's SLE process generates a family of random curves that grow in plane domains. The evolution is described by the classical Loewner differential equation with the driving function being $\sqrt\kappa B(t)$, where $B(t)$ is a standard Brownian motion and $\kappa$ is a positive parameter. SLE behaves differently for different value of $\kappa$. We use SLE$(\kappa)$ to emphasize the parameter. See \cite{LawSLE} and \cite{RS-basic} for the fundamental properties of SLE.

There are several versions of SLE, among which chordal SLE and radial SLE are most well known. They describe random curves that grow in simply connected domains. A number of statistical physics models in simply connected domains have been proved to converge in their scaling limits to chordal or radial SLE with different parameters.

People have been working on extending SLE to general plane domains. A version of SLE in doubly connected domains, called annulus SLE, was introduced in \cite{Zhan}. The definition uses annulus Loewner equation, in which the Poisson kernel function is used for the vector field, and the driving function is still $\sqrt\kappa B(t)$. Annulus SLE$(2)$ turns out to be the scaling limit of loop-erased random walk in doubly connected domains. In fact, loop-erased random walk in any finitely connected plane domain  converges to some SLE$(2)$-type curve (c.f.\ \cite{LERW}).

Annulus SLE defined in \cite{Zhan}  generates a trace in a doubly connected domain that starts from a marked boundary point and ends at a random point on the other boundary component (c.f.\ \cite{ann-prop}). This is different from the behavior of chordal SLE or radial SLE, whose trace ends at a fixed boundary point or interior point. The reason of this phenomena is that the definition of annulus SLE does not specify any point 
other than the initial point.

The annulus SLE$(\kappa;\Lambda)$ process was defined in \cite{whole} to describe SLE in doubly connected domains with one marked boundary point other than the initial point. Here the $\Lambda$ is a function, and the marked boundary point may or may not lie on the same boundary component as the initial point. The definition uses annulus Loewner equation with the driving function equal to $\sqrt\kappa B(t)$ plus some drift function. And the derivative of the drift function at any time is equal to the $\Lambda$ valued at the conformal type of the remaining domain together with the marked point and the tip of the SLE curve at that time.

There is very little restriction on the function $\Lambda$ in the above definition. For any $\kappa\in(0,4]$, there is a family of particular functions $\Lambda_{\kappa;\langle s\rangle}$, $s\in\R$, such that the annulus SLE$(\kappa;\Lambda_{\kappa;\langle s\rangle})$ process satisfies the remarkable reversibility properties as follows. Suppose $D$ is a doubly connected domain, and $z_0,w_0$ are two boundary points that lie on different boundary components. Let $\beta$ be an annulus SLE$(\kappa;\Lambda_{\kappa;\langle s\rangle})$ trace in $D$ that grows from $z_0$ with $w_0$ as the marked point. Then almost surely $\beta$ ends at $w_0$, and the time-reversal of $\beta$ is a time-change of an annulus SLE$(\kappa;\Lambda_{\kappa;\langle -s\rangle})$ trace in $D$ that grows from $w_0$ with $z_0$ as the marked point. This property was used (\cite{whole}) to prove the reversibility of whole-plane SLE$(\kappa)$ process for $\kappa\in(0,4]$.

In this paper we study the restriction property of the annulus SLE$(\kappa;\Lambda_{\kappa;\langle s\rangle})$ process.
We use $\mu_{\mbox{loop}}$ to denote the Brownian loop measure  defined in \cite{loop}, which is a $\sigma$-finite infinite measure on the space of loops, and define
 \BGE \cc=\cc(\kappa)=\frac{(6-\kappa)(3\kappa-8)}{2\kappa}.\label{cc}\EDE
 It is well known that $\cc$ is the central change for SLE$(\kappa)$. Set $\A_p=\{e^{-p}<|z|<1\}$, $\TT=\{|z|=1\}$ and $\TT_p=\{|z|=e^{-p}\}$. We will prove the following two theorems.

\begin{Theorem}
Let $p>0$, $\kappa\in(0,4]$, $s\in\R$, $z_0\in\TT$ and $w_0\in\TT_p$. Let $\nu $ be the distribution of an annulus SLE$(\kappa;\Lambda_{\kappa; \langle s\rangle})$ trace in $\A_p$ started from $z_0 $ with marked point $w_0 $. Let $L\subset\A_p$ be such that $\A_p\sem L$ is a doubly connected domain and $\dist(L,\{z_0,\TT_p\})>0$. Define a probability measure $\nu_{L}$ by
\BGE \frac{d\nu_{L }}{d\nu}=\frac {{\bf 1}_{\{\beta\cap L=\emptyset\}}}Z  \exp (\cc(\kappa) \mu_{\mbox{loop}}[{\cal L}_{L,p}]),\label{nu-L}\EDE
where $\beta$ is the SLE trace, 
$ {\cal L}_{L,p}$ is the set of all loops in $\A_p$ that intersect both $L$ and $\beta$, and $Z>0$ is a normalization factor. Then $\nu_{L }$ is the distribution of a time-change of  an annulus SLE$(\kappa;\Lambda_{\kappa; \langle s\rangle})$ trace in $\A_p\sem L$ started from $z_0$ with marked point $w_0$. \label{main-ann}
\end{Theorem}

\begin{Theorem}
 Let $p,\kappa,s,z_0,w_0,\nu$ be as in Theorem \ref{main-ann}. Let $L\subset\A_p$ be such that $\A_p\sem L$ is a simply connected domain, and $\dist(L,\{z_0,w_0\})>0$. Define $\nu_L$ by (\ref{nu-L}). Then $\nu_L$ is the distribution of a time-change of  a chordal SLE$(\kappa)$ trace in $\A_p\sem L$ from $z_0$ to $w_0$. \label{Main-strip}
\end{Theorem}

If $\kappa=\frac83$, then $\cc=0$. The above two theorems imply that, if we condition an annulus SLE$(\frac 83,\Lambda_{\frac83;\langle s\rangle})$ trace in $\A_p$ to avoid some set $L$, then the the resulting curve is a time-change of an annulus SLE$(\frac 83,\Lambda_{\frac83;\langle s\rangle})$ or chordal SLE$(\frac83)$ trace in $\A_p\sem L$. This is similar to the restriction property of chordal or radial SLE$(\frac8 3)$ (\cite{LSW-8/3}). If $\kappa\in(0,\frac83)$, then $\cc<0$, and the strong restriction property does not hold. But we may use the argument in \cite{LSW-8/3} to attach Brownian loops in $\A_p$ with density $-\cc$ to the trace to get a random shape with the restriction property.

The paper is organized as follows. We introduce notation, symbols and definitions in Section \ref{Sec2}, Section \ref{Sec3} and Section \ref{Sec4}. The proof of Theorem \ref{main-ann} is started at Section \ref{Sec5}, and finished at the end of Section \ref{Sec7}. The argument introduced in \cite{LSW-8/3} is used. In Section \ref{Sec8} we give a sketch of the proof of Theorem \ref{Main-strip}, and use Theorem \ref{Main-strip} to prove Theorem \ref{multiple}, which generates $n\ge 2$ mutually disjoint random curves crossing an annulus such that conditioned on all but one trace, the remaining trace is a chordal SLE$(\kappa)$ trace. 
We believe that, in the case $n=2$, if the inner circle of the annulus shrinks to a single point, then the two curves tend to the two arms of a two-sided radial SLE$(\kappa)$ (c.f.\ \cite{LawSLE}) in the disc. This may be used to understand the microscopic behavior of an SLE$(\kappa)$ trace near a typical point on this trace.

\section{Preliminary} \label{Sec2}
\subsection{Symbols and notation}
We will frequently use functions $\cot(z/2)$, $\tan(z/2)$, $\coth(z/2)$, $\tanh(z/2)$, $\sin(z/2)$, $\cos(z/2)$, $\sinh(z/2)$, and $\cosh(z/2)$. For simplicity, we write $2$ as a subscript. For example, $\cot_2(z)$ means $\cot(z/2)$, and $\cot_2'(z)=-\frac 12\sin_2^{-2}(z)$.

Let  $\TT=\{z\in\C:|z|=1\}$. For $p>0$, let $\A_p=\{z\in\C:1>|z|> e^{-p}\}$, $\St_p=\{z\in\C:0< \Imm z<p\}$, $\TT_p=\{z\in\C:|z|=e^{-p}\}$, and $\R_p=\{z\in\C:\Imm z=p\}$. Then $\pa\A_p=\TT\cup\TT_p$ and $\pa \St_p=\R\cup\R_p$. Let $e^i$ denote the map $z\mapsto e^{iz}$. Then $e^i$ is a covering map from $\St_p$ onto $\A_p$, maps $\R$ onto $\TT$ and maps $\R_p$ onto $\TT_p$.

A subset $K$ of a simply connected domain $D$ is called a hull in $D$ if $D\sem K$ is a simply connected domain.
A subset $K$ of a doubly connected domain $D$ is called a hull in $D$ if $D\sem K$ is a doubly connected domain, and $K$ is bounded away from a boundary component of $D$. In this case, we define
$\ccap_{D}(K):=\modd(D)-\modd(D\sem K)$
to be the capacity of $K$ in $D$, where $\modd(\cdot)$ is the modulus of a doubly connected domain. We have $0\le \ccap_{D}(K)<\modd(D)$, where the equality holds iff $K=\emptyset$. For example, the $L$ in Theorem \ref{main-ann} is a hull in $\A_p$.

We say a set $K\subset \C$ has period $p\in\C$ if $p+K=K$. We say that a function $f$ has progressive period $(p_1;p_2)$ if $f(\cdot\pm p_1)=f\pm p_2$. In this case, the definition domain of $f$ has period $p_1$, and the range of $f$ has period $p_2$.

An increasing function in this paper will always be strictly increasing.
For a real interval $J$, we use $C(J)$ to denote the space of real continuous functions on $J$.
The maximal solution to an ODE or SDE with initial value is the solution with the biggest definition domain.

A conformal map in this paper is an injective analytic function. We say that $f$ maps $D_1$ conformally onto $D_2$, and write $f:D_1\conf D_2$, if $f$ is a conformal map defined on the domain $D_1$ and $f(D_1)=D_2$. If, in addition, for $j=1,2$, $c_j$ is a point or a set in $D$ or on $\pa D$, and $f$ or its continuation maps $c_1$ onto $c_2$, then we write $f:(D_1;c_1)\conf(D_2;c_2)$.

Throughout this paper, a Brownian motion means a standard one-dimensional Brownian motion, and $B(t)$, $0\le t<\infty$, will always be used to denote a Brownian motion. This means that $B(t)$ is continuous, $B(0)=0$, and $B(t)$ has independent increment with $B(t)-B(s)\sim {\cal N}(0,t-s)$ for $t\ge s\ge 0$.

Many functions in this paper depend on two variables. The first variable represents time or modulus, and the second variable does not. We use $\pa_t$ and $\pa_t^{n}$ to denote the partial derivatives w.r.t.\ the first variable, and use $'$, $''$, and the superscripts $(h)$ to denote the partial derivatives w.r.t.\ the second variable.

\subsection{Special functions} \label{section-special}
 For $t>0$, define $$\SA(t,z)=\lim_{M\to\infty}\sum_{k=-M}^M \frac{e^{2kt}+z}{e^{2kt}-z}
=\PV\sum_{2\mid n} \frac{e^{nt}+z}{e^{nt}-z},$$ 
$$ \HA(t,z)=-i\SA(t,e^i(z))=-i\PV\sum_{2\mid n} \frac{e^{nt}+e^{iz}}{e^{nt}-e^{iz}}=\PV\sum_{2\mid n}\cot_2(z-int). $$ 
Then $\HA(t,\cdot)$ is a meromorphic function in $\C$, whose poles are $\{2m\pi+i2kt:m,k\in\Z\}$, which are
all simple poles with residue $2$. Moreover, $\HA(t,\cdot)$ is an odd function and takes real values on $\R\sem\{\mbox{poles}\}$;
$\Imm \HA(t,\cdot)\equiv -1$ on $\R_t$;  $\HA(t,\cdot)$ has period $2\pi$ and progressive period $(i2t;-2i)$.
Let $ \rA(t)\in\R$ 
be such that the power series expansion of $\HA(t,\cdot)$ near $0$ is
 \BGE \HA(t,z)=\frac 2z+ \rA(t)z+O(z^3),\label{Taylor}\EDE

Let $\SA_I(t,z)=\SA(t,e^{-t}z)-1$ and $\HA_I(t,z)=-i\SA_I(t,e^{iz})=\HA(t,z+it)+i$.
It is easy to check: $$ \SA_I(t,z)=\PV\sum_{2\nmid n} \frac{e^{nt}+z}{e^{nt}-z}, \quad
\HA_I(t,z) =\PV\sum_{2\nmid n} \cot_2(z-int). $$ 
So $\HA_I(t,\cdot)$ is a meromorphic function in $\C$ with poles $\{2m\pi+i(2k+1)t:m,k\in\Z\}$, which are
all simple poles with residue $2$;  $\HA_I(t,\cdot)$ is an odd function and takes real values on $\R$;  $\HA_I(t,\cdot)$ has period $2\pi$ and progressive period $(i2t;-2i)$.

 It is possible to express $\HA$ and $\HA_I$ using classical functions. Let $\theta(\nu,\tau)$ and $\theta_k(\nu,\tau)$, $k=1,2,3$, be the Jacobi theta functions defined in \cite{elliptic}. Define $\Theta(t,z)=\theta(\frac{z}{2\pi},\frac{it}\pi)$ and $\Theta_I(t,z)=\theta_2(\frac{z}{2\pi},\frac{it}\pi)$. Then $\Theta(t,\cdot)$ has antiperiod $2\pi$,  $\Theta_I(t,\cdot)$ has period $2\pi$, and
\BGE \HA=2\,\frac{\Theta'}{\Theta},\quad \HA_I=2\,\frac{\Theta_I'}{\Theta_I}. \label{HA-Theta}\EDE

It is useful to rescale the special functions. Let
\BGE  \ha \Theta(t,z)=e^{\frac{z^2}{4 t}}\Big(\frac\pi t\Big)^{\frac 12} \Theta\Big(\frac{\pi^2}t,\frac \pi tz\Big),\quad \ha \Theta_I(t,z)=e^{\frac{z^2}{4 t}}\Big(\frac\pi t\Big)^{\frac 12} \Theta_I\Big(\frac{\pi^2}t,\frac \pi tz\Big).\label{Theta-ha}\EDE
From the Jacobi identities, we have $\ha\Theta(t,z)=\theta(i\frac{z}{2\pi},\frac{it}\pi)=\Theta(t,iz)$ and $\ha\Theta_I(t,z)=\theta_1(i\frac{z}{2\pi},\frac{it}\pi)$. From the product representations of $\theta_1$, we get
\BGE \ha\Theta_I(t,z)=2e^{-\frac t4}\cosh_2(z)\prod_{m=1}^\infty (1-e^{-2mt})(1+e^{z-2mt})(1+e^{-z-2mt}).\label{Theta-ha-prod}\EDE

Let $\ha \HA=2\,\frac{\ha\Theta'}{\ha\Theta}$ and $\ha\HA_I=2\,\frac{\ha\Theta_I'}{\ha\Theta_I}$. 
From (\ref{HA-Theta}) and (\ref{Theta-ha}) we have
 \BGE \ha \HA(t,z)=\frac \pi t\HA\Big(\frac{\pi^2}t,\frac\pi t z\Big)+\frac zt,\qquad \ha\HA_I(t,z)=\frac \pi t\HA_I\Big(\frac{\pi^2}t,\frac\pi t z\Big)+\frac zt.\label{ha-HA}\EDE
Since $\ha\Theta(t,z)=\Theta(t,iz)$ and $\HA_I(t,z)=\HA(t,z+it)+i$, we have
\BGE \ha\HA(t,z)=i\HA(t,iz)=\PV\sum_{2\mid n}\coth_2(z-nt);\label{HA-coth}\EDE
\BGE \ha\HA_I(t,z)=\ha\HA(t,z+\pi i)=\PV\sum_{2\mid n}\tanh_2(z-nt).\label{ha-HA-exp}\EDE

From (\ref{HA-coth}), the power series expansion of $\ha\HA(t,\cdot)$ near $0$ is
 \BGE \ha\HA(t,z)=\frac 2z+ \ha\rA(t)z+O(z^3),\label{Taylor-ha}\EDE
where $ \ha\rA(t):=-\sum_{k=1}^\infty \sinh^{-2}(kt)+\frac 16=O(e^{-t})+\frac 16$ as $t\to\infty$. Hence we may define
\BGE \ha\RA(t)=-\int_t^\infty (\ha\rA(s)-\frac 16)ds,\quad 0<t<\infty.\label{RA-rA-ha}\EDE Then $\ha\RA$ is positive and decreasing as $\ha\rA-\frac 16<0$. From (\ref{Taylor}),  (\ref{ha-HA}), and (\ref{Taylor-ha}), we have
\BGE \ha\rA(t)=\Big(\frac \pi t\Big)^2\rA\Big(\frac {\pi^2} t\Big)+\frac 1t.\label{rA-ha}\EDE

 \section{Loewner equations} \label{Sec3}
 \subsection{Annulus Loewner equation} \label{section-ann-orig}
 The annulus Loewner equations are defined in \cite{Zhan}. Fix $p\in(0,\infty)$ and $T\in(0,p]$. Let $\xi\in C([0,T))$. The annulus Loewner equation of modulus $p$ driven by $\xi$ is
$$ \pa_t  g(t,z)= g(t,z) \SA(p-t, g(t,z)/e^{i\xi(t)}),\quad  g(0,z)=z.$$  
For $0\le t<T$, let $K(t)$ denote the set of $z\in\A_p$ such that the
solution $ g(s,z)$ blows up before or at time $t$. Then each $K(t)$ is a hull in $\A_p$, $\ccap_{\A_p}(K(t))=t$, and $g(t,\cdot)$ maps $\A_p\sem K(t)$ conformally onto $\A_{p-t}$, and maps $\TT_p$ onto $\TT_{p-t}$. We call $K(t)$ and $ g(t,\cdot)$, $0\le t<T$, the annulus
Loewner hulls and maps of modulus $p$ driven by $\xi$.

It is known that, if $\xi$ is a semi-martingale whose stochastic part is $\sqrt\kappa B(t)$, and whose drift part is continuously differentiable, then
$\xi$ generates an annulus Loewner trace $\beta$ of modulus $p$, which means that \BGE \beta(t):=\lim_{\A_{p-t}\ni z\to e^{i\xi(t)}} g(t,\cdot)^{-1}(z)\label{trace}\EDE
exists for all $0\le t<T$, and $\beta$ is a continuous simple curve in $\A_p\cup\TT$ with $\beta(0)=e^{i\xi(0)}\in\TT$. If $\kappa\in(0,4]$, then $\beta$ is simple and $\beta((0,T))\subset \A_p$. In this case, $K(t)=\beta((0,t])$ for $0\le t<T$, and we say that $\beta$ is parameterized by its capacity in $\A_p$ w.r.t.\ $\TT_p$, i.e., $\ccap_{\A_p}(\beta((0,t]))=t$ for $0\le t<T$. 

On the other hand, if $\beta(t)$, $0\le t<T$, is a simple curve with $\beta(0)\in\TT$, $\beta((0,T))\subset\A_p$, and if $\beta$ is parameterized by its capacity in $\A_p$ w.r.t.\ $\TT_p$, then $\beta$ is a simple annulus Loewner trace of modulus $p$ driven by some $\xi\in C([0,T))$. If $\beta$ is not parameterized by its capacity, then $\beta(v^{-1}(t))$, $0\le t<v(T)$, is an annulus Loewner trace of modulus $p$, where $v(t):=\ccap_{\A_p}(\beta((0,t]))$ is an increasing function with $v(0)=0$.

\subsection{Covering annulus Loewner equation}
The covering annulus Loewner equation of modulus $p$ driven by $\xi\in C([0,T))$ is
\BGE \pa_t{\til g}(t,z)=\HA(p-t,\til g(t,z)-\xi(t)),\quad \til g(0,z)=z.\label{annulus-eqn-covering*}\EDE
For $0\le t<T$, let $\til K(t)$ denote the set of $z\in\St_p$ such that the
solution $\til g(s,z)$ blows up before or at time $t$.
Then for $0\le t<T$,
\BGE \til g(t,\cdot):(\St_p\sem \til K(t);\R_p)\conf (\St_{p-t};\R_{p-t}).\label{conf-g-til}\EDE
We call $\til K(t)$ and $\til  g(t,\cdot)$, $0\le t<T$, the covering annulus
Loewner hulls and maps of modulus $p$ driven by $\xi$.

The relation between the covering annulus Loewner equation and the annulus Loewner equation is as follows.
Let $K(t)$ and $ g(t,\cdot)$ be the annulus Loewner hulls and maps of modulus $p$ driven by $\xi$.
 Then we have $\til K(t)=(e^i)^{-1}(K(t))$ and  $e^i\circ \til g(t,\cdot)= g(t,\cdot)\circ e^i$,  $0\le t<T$. 
Thus, $\til K(t)$ has period $2\pi$, and $\til g(t,\cdot)$ has progressive period $(2\pi;2\pi)$.

If $\xi$ generates an annulus Loewner trace $\beta$ defined by (\ref{trace}), then there is a continuous simple curve $\til\beta(t)$, $0\le t<T$,
 which is defined by \BGE \til\beta(t)=\lim_{\St_{p-t}\ni z\to{\xi(t)}} \til g(t,\cdot)^{-1}(z),\quad0\le t<T.\label{trace-til*}\EDE
Such $\til\beta$ is called the covering annulus Loewner trace of modulus $p$ driven by $\xi$, and satisfies that
 $\beta=e^i\circ \til \beta$ and $\til\beta(0)=\xi(0)$. If $\beta$ is simple with $\beta((0,T))\subset\A_p$, then $\til\beta$ is also simple, $\til\beta((0,T))\subset\St_p$, and $\til K(t)=\til\beta((0,t])+2\pi\Z$, $0\le t<T$.

Since $\til g(t,\cdot)$ maps $\R_p$ onto $\R_{p-t}$ and $\HA_I(t,z)=\HA(t,z+it)+i$, we have
\BGE \pa_t \Ree \til g(t,z)=\HA_I(p-t,\Ree\til g(t,z)-\xi(t)),\quad z\in\R_p.\label{ODE-HA-I*}\EDE
Differentiating  (\ref{ODE-HA-I*}) w.r.t.\ $z$, we see that
\BGE \pa_t {\til g}'(t,z)=\til g'(t,z)\HA_I'(p-t,\Ree\til g(t,z)-\xi(t)),\quad z\in\R_p.\label{deriv2*}\EDE

Since $\SA(p-t,\cdot)$ and $\HA(p-t,\cdot)$ have period $2\pi$, for any $n\in\Z$, $\xi$ and $\xi+2n\pi$ generate  the same family of annulus Loewner maps and the same family of covering annulus Loewner maps.

\subsection{Strip Loewner evolution}
Strip Loewner equations will be used in Section \ref{Sec8}. The strip Loewner equation (\cite{thesis}) driven by $\xi\in C([0,T))$ is
$$\pa_t \til g(t,z)=\coth_2(\til g(t,z)-\xi(t)),\quad 0\le t<T,\quad \til g(0,z)=z.$$
For $0\le t<T$, let $\til K(t)$ denote the set of $z\in\St_\pi$ such that the
solution $ \til g(s,z)$ blows up before or at time $t$. Then $\til K(t)$ and $\til g(t,\cdot)$, $0\le t<T$, are called the strip
Loewner hulls and maps driven by $\xi$. For each $t\in[0,T)$, $\til K(t)$ is a bounded hull in $\R_\pi$ with $\dist(\til K(t),\R_\pi)>0$, $\til g(t,\cdot):(\St_\pi\sem \til K(t);\R_\pi)\conf (\St_{\pi};\R_\pi)$, and $\til g(t,z)-z\to \pm t$ as $z\to\pm\infty$ in $ \St_\pi\sem \til K(t)$.
If $\til K$ is a bounded hull in $\R_\pi$ with $\dist(\til K(t),\R_\pi)>0$, then there exist a number $c_{\til K}\ge 0$ and a map $\til g_{\til K}$ determined by $\til K$ such that $\til g_{\til K}:(\St_\pi\sem \til K;\R_\pi)\conf (\St_{\pi};\R_\pi)$ and  $\til g_{\til K}-z\to \pm c_{\til K}$ as $z\to\pm\infty$. We call $c_{\til K}$ the capacity of $\til K$ in $\St_\pi$ w.r.t.\ $\R_\pi$. Thus, the capacity of $\til K(t)$ in $\St_\pi$ w.r.t.\ $\R_\pi$ is $t$, and $\til g(t,\cdot)=\til g_{\til K(t)}$.

Since $\til g(t,\cdot)$ maps $\R_\pi$ onto $\R_\pi$ and $\coth_2(z+\pi i)=\tanh_2(t,z)$, we have
\BGE \pa_t \Ree \til g(t,z)=\tanh_2(\Ree\til g(t,z)-\xi(t)),\quad z\in\R_\pi.\label{ODE-HA-I-strip}\EDE
Differentiating  (\ref{ODE-HA-I-strip}) w.r.t.\ $z$, we see that
\BGE \pa_t {\til g}'(t,z)=\til g'(t,z)\tanh_2'(\Ree\til g(t,z)-\xi(t)),\quad z\in\R_\pi.\label{deriv2-strip}\EDE

If $\xi$ is a semi-martingale whose stochastic part is $\sqrt\kappa B(t)$, and whose drift part is continuously differentiable, then
$\xi$ generates a strip Loewner trace $\til\beta$, which is defined by
\BGE \til \beta(t):=\lim_{\St_\pi\ni z\to \xi(t)} \til g(t,\cdot)^{-1}(z), \quad 0\le t<T.  \label{trace-strip}\EDE
Such $\til\beta$ is a continuous curve in $\St_\pi\cup\R$ which satisfies that $\til\beta(0)=\xi(0)\in\R$. If $\kappa\in(0,4]$, then $\til\beta$ is simple, $\til\beta((0,T))\subset\St_\pi$, and $\til K(t)=\til\beta((0,t])$ for $0\le t<T$.

 On the other hand, suppose $\til\beta(t)$ is a simple curve in $\St_\pi\sem\R$, which intersects $\R$ only at $t=0$. Let $v(t)$ be the capacity of $\til\beta((0,t])$ in $\St_\pi$ w.r.t.\ $\R_\pi$. Then $v$ is a continuous increasing function, which maps $[0,T)$ onto $[0,S)$ for some $S\in(0,\infty]$, and there is $\xi\in C([0,S))$ which generates the  strip Loewner trace $\til\beta\circ v^{-1}$.

The chordal SLE$(\kappa;\rho)$ process defined in \cite{LSW-8/3} naturally extends to strip SLE$(\kappa;\rho)$ process. Let $\kappa>0$ and $\rho\in\R$. Let $x_0,y_0\in\R$. Let $\xi(t)$ and $q(t)$, $0\le t<\infty$, be the solution of
$$d\xi(t)=\sqrt\kappa dB(t)+\frac\rho 2\tanh_2(\xi(t)-q(t))dt,\quad \xi(0)=x_0;$$
$$dq(t)=\tanh_2(q(t)-\xi(t)),\quad q(0)=y_0.$$
Then the strip Loewner trace $\til\beta$ driven by $\xi$ is called a strip SLE$(\kappa;\rho)$ trace in $\St_\pi$ started from $x_0$ with marked point $y_0+\pi i$. 
From \cite{SW} we know that, when $\rho=\kappa-6$, $\til\beta$ is a time-change of a chordal SLE$(\kappa)$ trace in $\St_\pi$ from $x_0$ to $y_0+\pi i$, stopped when it hits $\R_\pi$. If, in addition, $\kappa\le 4$, since the chordal SLE$(\kappa)$ trace does not hit $\R_\pi$ before it ends, we see that $\til\beta$ is a time-change of a complete chordal SLE$(\kappa)$ trace.

\section{One SLE Curve Crossing an Annulus}\label{Sec4}
\subsection{Annulus SLE with one marked point}\label{section-chordal}
We now cite some definitions in Section 4.1 of \cite{whole}.
\begin{Definition} A covering crossing annulus drift function is a real valued $C^{0,1}$ differentiable function defined on  $(0,\infty)\times\R$. A covering crossing annulus drift function with period $2\pi$ in its second variable is called a crossing annulus drift function.
\label{drift-function}
\end{Definition}

\begin{Definition}
Suppose $\Lambda$ is a covering crossing annulus drift function. Let $\kappa> 0$, $p>0$, and $x_0,y_0\in\R$. Let $\xi(t)$, $0\le t< p$, be the maximal solution to the SDE
\BGE d \xi(t)= \sqrt\kappa dB(t)+\Lambda(p-t,\xi(t)-\Ree \til g(t,y_0+ pi))dt,\quad  \xi(0)=x_0,\label{xi-crossing}\EDE
where $\til g(t,\cdot)$, $0\le t<p$, are the covering annulus Loewner maps of modulus $p$ driven by $\xi$. Then the covering annulus Loewner trace of modulus $p$ driven by $\xi$ is called the covering  annulus SLE$(\kappa;\Lambda)$ trace in $\St_p$ started from $x_0$ with marked point $y_0+p i$. \label{covering-ann-crossing}
\end{Definition}


\begin{Definition} Suppose $\Lambda$ is a crossing annulus drift function. Let $\kappa\ge 0$, $p>0$, $a\in\TT$
and $b\in\TT_p$. Choose  $x_0,y_0\in\R$ such that $a=e^{ix_0}$ and $b=e^{-p+iy_0}$.   Let $\xi(t)$, $0\le t< p$, be the maximal solution to (\ref{xi-crossing}). The  annulus Loewner trace of modulus $p$ driven by $\xi(t)$, $0\le t<p$, is called the annulus SLE$(\kappa;\Lambda)$ trace in $\A_p$ started from $a$ with marked point $b$.
 \label{crossing-ann}
\end{Definition}

\no {\bf Remark.} The above definition does not depend on the choices of $x_0$ and $y_0$ because $\Lambda(p-t,\cdot)$ has period $2\pi$, $\til g(t,\cdot) $ has progressive period $(2\pi;2\pi)$, and for any $n\in\Z$, the annulus Loewner objects driven by $\xi(t)+2n\pi$ agree with those driven by $\xi(t)$. Via conformal maps,   we can define annulus SLE$(\kappa;\Lambda)$ trace in any doubly connected domain. 

\subsection{Annulus SLE with reversibility} \label{section-limit}
A family of functions are defined in Section 7 of \cite{whole}, which are $\ha\Psi_\infty$, $\ha\Psi_q$, $\ha\Psi_0$, $\Psi_0$, $\Psi_m$, $m\in\Z$, $\Psi_{\langle s\rangle}$, $\Lambda_0$, and $\Lambda_{\langle s\rangle}$, $s\in\R$. They are all smooth functions on $(0,\infty)\times\R$, and depend on three parameters: $\kappa\in(0,4]$, $\sigma\in [0,\frac 4\kappa)$, and $\tau=\frac\kappa4-\sqrt{\frac{\kappa^2}{16}+\kappa\sigma}\le 0$. Now we suppose $\kappa\in(0,4]$ is fixed, and
\BGE \sigma=\frac 4\kappa-1\ge 0,\quad \tau=\frac\kappa 2-2\le 0.\label{sigma-tau}\EDE
Then these function depend only on $\kappa\in(0,4]$, $m\in\Z$ and $s\in\R$. For simplicity, we omit the symbol $\kappa$.
 The $\Lambda_{\langle s\rangle}$ here is the $\Lambda_{\kappa;\langle s\rangle}$ in Theorem \ref{main-ann} and Theorem \ref{Main-strip}.

The $\ha\Psi_\infty$ is defined in (7.31) of \cite{whole}:
\BGE \ha\Psi_\infty(t,x)=e^{-\frac{\tau^2t}{2\kappa}}\cosh_2^{\frac 2\kappa\tau}(x).\label{psiqinfty}\EDE
The $\ha \Psi_q$ is defined by (7.33) of \cite{whole}:
\BGE \ha \Psi_q(t,x)=\EE\Big[\exp\Big(\sigma\int_0^\infty \ha\HA_{I,q}'(t+s,X_x(s))ds\Big)\Big],\label{psiq}\EDE
where $\ha\HA_{I,q}$ is defined by (7.8) of \cite{whole}:
$ \ha\HA_{I,q}(t,z)=\ha\HA_I(t,z)-\tanh_2(z)$, 
and $X_x(t)$, $0\le t<\infty$, is a diffusion process which satisfies SDE (7.2) of \cite{whole}:
\BGE dX_x(t)=\sqrt\kappa dB(t)+\tau \tanh_2(X_x(t))dt,\quad X_x(0)=x.\label{dX}\EDE
The $\ha\Psi_0$ is defined in Theorem 7.2 of \cite{whole}:
\BGE \ha\Psi_0=\ha\Psi_\infty\ha\Psi_q.\label{ha-Psi-0}\EDE
The $\Psi_0$ is defined in Theorem 7.3 of \cite{whole}:
\BGE \Psi_0(t,x)=e^{-\frac{x^2}{2\kappa t}}\Big(\frac\pi t\Big)^{\sigma+\frac 12}\ha\Psi_0\Big(\frac{\pi^2}t,\frac \pi tx\Big).\label{Psi-ha}\EDE
For $m\in\Z$ and $s\in\R$, the $\Psi_m$ and $\Psi_{\langle s\rangle}$ are defined in Theorem 7.4 of \cite{whole}:
\BGE \Psi_m(t,x)=\Psi_0(t,x-2m\pi),\quad \Psi_{\langle s\rangle} =\sum_{m\in\Z} e^{\frac{2\pi}\kappa ms}\Psi_m.\label{Psi-langle}\EDE
The functions $\ha\Psi_\infty$, $\ha\Psi_q$, $\ha\Psi_0$, $\Psi_0$, $\Psi_m$, $\Psi_{\langle s\rangle}$ are all positive. The functions $\Lambda_0$ and $\Lambda_{\langle s\rangle}$ are defined in Proposition 7.4 and Theorem 7.4, respectively, of \cite{whole}: $ \Lambda_0=\kappa\frac{\Psi_0'}{\Psi_0}-\HA_I$,  $\Lambda_{\langle s\rangle} =\kappa \frac{\Psi_{\langle s\rangle}'}{\Psi_{\langle s\rangle}}-\HA_I$. 
 For the sake of completeness, we now  define $\Lambda_m=\kappa\frac{\Psi_m'}{\Psi_m}-\HA_I=\Lambda_0(\cdot-2m\pi)$ and
\BGE \Gamma_m=\Psi_m\Theta_I^{-\frac 2\kappa},\quad \Gamma_{\langle s\rangle}=\Psi_{\langle s\rangle}\Theta_I^{-\frac 2\kappa}.\label{Gamma-Psi}\EDE
From (\ref{HA-Theta}), we see that $\Lambda_m$ and $\Lambda_{\langle s\rangle}$ have simpler expressions:
\BGE \Lambda_m=\kappa\frac{\Gamma_m'}{\Gamma_m},\quad \Lambda_{\langle s\rangle}=\kappa\frac{\Gamma_{\langle s\rangle}'}{\Gamma_{\langle s\rangle}}.\label{Lambda-Gamma}\EDE
From Lemma 5.2 of \cite{whole}, we see that $\Gamma_m$ and $\Gamma_{\langle s\rangle}$ solve the PDE (5.6) in \cite{whole}. Since we here set the value of $\sigma$ by (\ref{sigma-tau}), this PDE becomes (5.2) in \cite{whole}, i.e.,
 \BGE \pa_t\Gamma_m=\frac\kappa 2 \Gamma_m''+\Gamma_m'\HA_I+\alpha \HA_I'\Gamma_m, \label{PDE-Gamma}\EDE
where \BGE \alpha=\frac{6-\kappa}{2\kappa}.\label{alpha}\EDE

Define $\ha\Gamma_0$ on $(0,\infty)\times\R$ such that
\BGE  \ha\Gamma_0(t,x)=\Big(\frac \pi t\Big)^{\alpha}\Gamma_0\Big(\frac{\pi^2}t,\frac\pi t x\Big).\label{ha-Gamma}\EDE
From (\ref{Theta-ha}), (\ref{Psi-ha}), and (\ref{Gamma-Psi}), we have
\BGE \ha\Gamma_0=\ha\Psi_0\ha\Theta_I^{-\frac 2\kappa}.\label{Gamma-m-ha}\EDE
Define $\ha\Theta_{I,\infty}$, $\ha\Theta_{I,q}$, $\ha\Gamma_\infty$, and $\ha\Gamma_q$ on $(0,\infty)\times\R$ such that
\BGE \ha\Theta_{I,\infty}(t,x)=2e^{-\frac t4}\cosh_2(x);\quad  \ha\Theta_{I,q}=\ha\Theta_I/\ha\Theta_{I,\infty}; \label{ha-Theta-infty}\EDE
\BGE \ha\Gamma_{\infty}(t,x)= 2^{-\frac 2\kappa} e^{-\frac{\tau^2-1}{2\kappa}t}\cosh_2(x)^{\frac 2\kappa(\tau-1)};\quad \ha\Gamma_q=\ha\Gamma_0/\ha\Gamma_\infty. \label{ha-Gamma-infty}\EDE
One may check that $\ha\Gamma_\infty$ solves
\BGE -\pa_t \ha\Gamma_\infty=  \frac\kappa 2 \ha\Gamma_\infty''+\ha\Gamma_\infty'\tanh_2+\alpha \tanh_2'\ha\Gamma_\infty.\label{PDE-Gamma-infty}\EDE
From (\ref{psiqinfty}) we have $\ha\Gamma_\infty=\ha\Psi_\infty\ha\Theta_{I,\infty}^{-\frac 2\kappa}$.
From (\ref{ha-Psi-0}) and (\ref{Gamma-Psi}) we have
\BGE \ha\Gamma_q=\ha\Psi_q\ha\Theta_{I,q}^{-\frac 2\kappa}.\label{Gamma-q}\EDE

Let $p>0$ and $x_0,y_0\in\R$. Let $y_m=y_0+2m\pi$, $m\in\Z$. Consider the following two SDEs.
\BGE d\xi(t)=\sqrt\kappa dB(t)+\Lambda_0(p-t,\xi(t)-\Ree \til g(t,y_m+p i))dt,\quad 0\le t<p,\quad \xi(0)=x_0,\label{xi-p-m}\EDE
\BGE d\xi (t)=\sqrt\kappa dB(t)+\Lambda_{\langle s\rangle}(p-t,\xi(t)-\Ree \til g (t,y_0+p i))dt,\quad 0\le t<p,\quad\xi(0)=x_0,\label{xi-p}\EDE
where $\til g (t,\cdot)$  are the covering annulus Loewner maps driven by $\xi $. Let $\mu_m$ or $\mu_{\langle s\rangle}$ denote the distribution of $(\xi(t),0\le t<p)$ if it solves (\ref{xi-p-m}) or (\ref{xi-p}), respectively. Then
\BGE \mu_{\langle s\rangle}=\sum_{m\in\Z}e^{\frac{2\pi}\kappa ms}\, \frac{ \Psi_m(p,x_0-y_0)}{\Psi_{\langle s\rangle}(p,x_0-y_0)}\,\mu_m=\sum_{m\in\Z}e^{\frac{2\pi}\kappa ms}\, \frac{ \Gamma_0(p,x_0-y_m)}{\Gamma_{\langle s\rangle}(p,x_0-y_0)}\, \mu_m,\label{convex}\EDE
where the first equality follows from Proposition 7.4 in \cite{whole}, and the second equality follows from (\ref{Psi-langle}), (\ref{Gamma-Psi}), and the fact that $\Theta_I(p,\cdot)$ has period $2\pi$.

Let $\beta$ and $\til\beta$ be the annulus Loewner trace and covering annulus Loewner trace, respectively, of modulus $p$, driven by $\xi$. If $(\xi)$ has distribution $\mu_m$, then  $\til\beta$ is a covering annulus SLE$(\kappa;\Lambda_0)$ trace in $\St_p$ started from $x_0$ with marked point $y_m+p i$. If $(\xi)$ has distribution $\mu_{\langle s\rangle}$, then  $\beta$ is an annulus SLE$(\kappa;\Lambda_{\langle s\rangle})$ trace in $\A_p$ started from $e^{ix_0}$ with marked point $e^{iy_0-p}$. Let ${\cal E}_m$ denote the event that the covering trace ends at $y_m + p i$.
Proposition 7.4, Theorem 8.3, and Theorem 9.3 in \cite{whole}  together imply that  $\mu_m({\cal E}_m)=1$ and $\mu_{\langle s\rangle}(\bigcup_{m\in\Z}{\cal E}_m)=1$. Since ${\cal E}_m$, $m\in\Z$, are mutually disjoint, the $\mu_m$'s are singular to each other. From (\ref{convex}) we have
\BGE \frac{d\mu_m}{d\mu_{\langle s\rangle}}=e^{\frac{2\pi}\kappa ms}\, \frac{ \Gamma_0(p,x_0-y_m)}{\Gamma_{\langle s\rangle}(p,x_0-y_0)}\,{\bf 1}_{{\cal E}_m}.\label{dmu-m/dmu}\EDE


\subsection{Some estimations}
\begin{Lemma}
  For any   $t>0$ and $0\le x\le 3t$,
  $$\ha\HA_{I,q}'(t,x)<
\min\Big\{\frac 12,2e^{x-2t}\Big\}+\frac{4e^{-t}}{1-e^{-2t}}.$$
  \label{haIq<}
\end{Lemma}
{\bf Proof.} Since $\ha\HA_{I,q}(t,z)=\ha\HA_I(t,z)-\tanh_2(z)$,  from (\ref{ha-HA-exp}), we have $$\ha\HA_{I,q}'(t,x)=\tanh_2'(x-2t)+\sum_{n=2}^\infty \tanh_2'(x-2nt)+\sum_{n=-\infty}^{-1} \tanh_2'(x-2nt).$$
Note that $ \tanh_2'(x)=\frac{2}{(e^{x/2}+e^{-x/2})^2}\le  \min\Big\{\frac 12,2e^{x},2e^{-x}\Big\}$ for $x\in\R$.
If $0\le x\le 3t$, then
$$\sum_{n=2}^\infty \tanh_2'(x-2nt)\le 2\sum_{n=2}^\infty e^{x-2nt}=\frac{2e^{x-4t}}{1-e^{-2t}}\le \frac{2e^{-t}}{1-e^{-2t}},$$
$$\sum_{n=-\infty}^{-1} \tanh_2'(x-2nt)\le 2\sum_{n=-\infty}^{-1} e^{2nt-x}=\frac{2e^{-x-2t}}{1-e^{-2t}}\le \frac{2e^{-2t}}{1-e^{-2t}}.$$
The conclusion follows from the above displayed formulas. $\Box$

\begin{Proposition} If $F$ is one of the following functions: $\ha\Theta_{I,q}$, $\ha\Psi_q$, or $\ha\Gamma_q$, then
\begin{enumerate}
  \item [(i)] $\lim_{2t-|x|\to+\infty} \ln(F(t,x))=0$;
    \item [(ii)] for every $R>0$, $\ln(F)$ is bounded on $\{t\ge R,|x|\le 2t+R\}$;
\end{enumerate} \label{est-Gamma}
\end{Proposition}
{\bf Proof.}   From (\ref{Theta-ha-prod}) and (\ref{ha-Theta-infty}), the conclusion is clearly true for $F=\ha\Theta_{I,q}$. From (\ref{Gamma-q}), we suffice to prove this proposition for $F=\ha\Psi_q$. Throughout this proof, we use $O_t(1)$ to denote a positive quantity which depends on $\kappa,\sigma,t$, and is uniformly bounded when $t$ is bigger than any positive constant.

Fix $t>0$ and $x\in\R$. Let $X_x(s)$ be as in (\ref{dX}), and $(\F_s)$ be the filtration generated by $(X_x(s))$.
Define a uniformly integrable martingale $M_{t,x}(s)$, $0\le s<\infty$, by
$$ M_{t,x}(s):=\EE\Big[\exp\Big(\sigma\int_0^\infty \ha\HA_{I,q}'(t+r,X_x(r))dr\Big)\Big|\F_s\Big].$$
From (\ref{psiq}) we have $M_{t,x}(s)=\ha\Psi_q(t+s,X_x(s))\exp\Big(\sigma\int_0^s \ha\HA_{I,q}'(t+r,X_x(r))dr\Big)$. 
Suppose $S$ is an a.s.\ finite $(\F_s)$-stopping time. From the Optional Stopping Theorem, $\EE[M_{t,x}(S)]=M_{t,x}(0)$. Since $M_{t,x}(0)=\ha\Psi_q(t,x)$, we have
\BGE \ha\Psi_q(t,x)=\EE\Big[\ha\Psi_q(t+S,X_x(S))\exp\Big(\sigma\int_0^S \ha\HA_{I,q}'(t+s,X_x(s))ds\Big)\Big].\label{opt}\EDE

Let  $\lambda(s,x)=\min\{\frac 12,2e^{x-2s}\}$. If $0\le X_x(s)\le 3(t+s)$ for $0\le s\le S$, then from Lemma \ref{haIq<}, we have
$$\int_0^S \ha\HA_{I,q}'(t+s,X_x(s))ds\le \int_0^S \lambda(t+s,X_x(s))ds +\int_0^\infty \frac{4e^{-(t+s)}}{1-e^{-2(t+s)}}ds$$
$$=  \int_0^S \lambda(t+s,X_x(s))ds + 2\ln\Big(\frac{1+e^{-t}}{1-e^{-t}}\Big),$$
which together with (\ref{opt}) implies that
\BGE  \ha\Psi_q(t,x)\le \exp({O_t(1)e^{-t}}) \EE\Big[\ha\Psi_q(t+S,X_x(S))\exp\Big(\sigma\int_0^S \lambda(t+s,X_x(s))ds\Big)\Big].\label{opt2}\EDE

Recall that $\sigma\in[0,\frac4\kappa)$. Let $\sigma'=\frac\kappa 4\sigma$. From Proposition 7.1 in \cite{whole}, for any $c_0\in(1+\sigma',2)$, there is $C>0$ depending only on $\kappa$, $\sigma$, and $c_0$ such that for any $t\in(0,\infty)$ and $x\in\R$,
\BGE 1\le \ha\Psi_q(t,x)\le \exp\Big(C(t^{-1}+1) e^{(c_0-2)t}\Big)(1+Ce^{\frac 2\kappa|x|-\frac 2\kappa c_0t}).\label{psiq1}\EDE
This immediately implies that $\ln (\ha\Psi_q)$ is bounded on $\{|x|\le c_0 t,t\ge t_0\}$ for any $t_0>0$.

Choose any $c_0\in(1+\sigma',2)$ such that $c_0\ge\frac{2}{1+2/\kappa}$ and $c_0\ne 3-\frac \kappa 2$. Let $a=3-c_0\in (1,2-\sigma')$. Then $a\ne \frac\kappa 2$. Since $c_0-1>\sigma'$ and $2>a>1$, we have $a(2-a)=a(c_0-1)>\sigma'$. Thus,
\BGE -\frac 2\kappa a+\frac 2\kappa\frac{a+\sigma'}{c_0}=-\frac 2{\kappa c_0} (a(c_0-1)-\sigma')<0;\label{ac01}\EDE
\BGE -\frac 2\kappa a+\frac{\sigma}{2(2-a)}=-\frac 2{\kappa (2-a)} (a(2-a)-\sigma')<0 .\label{ac02}\EDE
For $m\in\N\cup\{0\}$, let ${\cal G}_m$ denote the event that $\sqrt \kappa B(s)< a s+m$ for any $s\ge 0$.
Then $\emptyset={\cal G}_0\subset {\cal G}_1\subset\cdots\subset{\cal G}_m\subset{\cal G}_{m+1}\subset\cdots$. It is well known that $\PP[\bigcup_{m=0}^\infty {\cal G}_m]=1$ and
\BGE \PP[{\cal G}_m^c]\le e^{-\frac 2\kappa ma},\quad m\in\N.\label{PE<*}\EDE


Suppose $t>0$ and $2t\le x\le 3t$. Let $S$ be the first time that $X_x(s)\le 0$ or $X_x(s)\ge 3(t+s)$. Then $S$ is a stopping time, and $0\le X_x(s)\le 3(t+s)$ for $0\le s\le S$. Since $X_x(s)$ is recurrent, $S$ is a.s.\ finite. Since $\tau\le 0$ and $\tanh_2(x)\ge 0$ for $x\ge 0$, from (\ref{dX}) we have
\BGE X_x(s)\le x+as+m,\quad 0\le s\le S,\quad\mbox{on}\quad{\cal G}_m.\label{Em<*}\EDE
Let ${\cal E}_l$ and ${\cal E}_r$ denote the event that $X_x(S)=0$ and $X_s(S)=3(t+S)$, respectively. From (\ref{psiq1}) and the facts that $0>c_0-2\ge -\frac 2\kappa c_0$ and $3-c_0>0$ we see that
 \BGE \ha\Psi_q(t+S,X_x(S)) \le  \exp(O_t(1) e^{(c_0-2)t})\le O_t(1)\quad\mbox{on}\quad{\cal E}_l,\label{El}\EDE
\BGE \ha\Psi_q(t+S,X_x(S)) \le O_t(1) e^{\frac2\kappa(3-c_0)(t+S)}\quad\mbox{on}\quad{\cal E}_r.\label{Er}\EDE

From (\ref{opt2}) we have
$$\ha\Psi_q(t,x)\le  \exp(O_t(1)e^{-t})\sum_{m=1}^{\infty} \EE\Big[{\bf 1}_{({\cal G}_m\sem {\cal G}_{m-1})\cap{\cal E}_r}\ha\Psi_q(t+S,X_x(S))\exp\Big(\sigma\int_0^S \lambda(t+s,X_x(s))ds\Big)\Big]$$ \BGE + \exp(O_t(1)e^{-t})\sum_{m=1}^{\infty} \Big[{\bf 1}_{({\cal G}_m\sem {\cal G}_{m-1})\cap{\cal E}_l}\ha\Psi_q(t+S,X_x(S))\exp\Big(\sigma\int_0^S \lambda(t+s,X_x(s))ds\Big)\Big] \label{CE*}\EDE
Suppose ${\cal G}_m\cap{\cal E}_r$ occurs. From (\ref{Em<*}) we have $3(t+S)=X_x(S)\le x+aS+m$. Since $3-a=c_0>0$, we have
\BGE S\le \frac{x-3t+m}{c_0}\quad \mbox{on}\quad {\cal G}_m\cap{\cal E}_r.\label{S<}\EDE
Since $S\ge 0$, we see that ${\cal G}_m\cap{\cal E}_r=\emptyset$ when $m<3t-x$. Let $m_0=\lceil 3t-x\rceil$. Then
From (\ref{PE<*}), (\ref{Er}), and the fact that $\lambda\le \frac12$, we find that for any $m\in\N$ and $m\ge m_0$,
$$\EE\Big[{\bf 1}_{({\cal G}_m\sem {\cal G}_{m-1})\cap{\cal E}_r}\ha\Psi_q(t+S,X_x(S))\exp\Big(\sigma\int_0^S \lambda(t+s,X_x(s))ds\Big)\Big]$$
$$\le \EE\Big[e^{-\frac2\kappa(m-1)a} O_t(1) \exp\Big(\frac2\kappa(3-c_0)(t+S)+\frac\sigma 2 S\Big)\Big]$$
Since $\frac 2\kappa (3-c_0)+\frac\sigma 2=\frac 2\kappa(a+\sigma')>0$ , from (\ref{S<}) we find that the RHS of the above formula is
$$\le O_t(1)  \exp\Big(-\frac2\kappa a(m-1-t) +\frac 2\kappa(a+\sigma') \cdot \frac{x-3t+m}{c_0}\Big) $$
So we have
$$ \sum_{m=1}^{\infty} \EE\Big[{\bf 1}_{({\cal G}_m\sem {\cal G}_{m-1})\cap{\cal E}_r}\ha\Psi_q(t+S,X_x(S))\exp\Big(\sigma\int_0^S \lambda(t+s,X_x(s))ds\Big)\Big]$$
$$\le O_t(1)\sum_{m=m_0}^\infty  \exp\Big(-\frac2\kappa a(m-t) +\frac 2\kappa(a+\sigma') \cdot \frac{x-3t+m}{c_0}\Big) $$
$$= O_t(1)\exp\Big(\frac2\kappa at+\frac 2\kappa\frac{a+\sigma'}{c_0}(x-3t)\Big)\sum_{m=m_0}^\infty \exp\Big(-\frac 2\kappa a+\frac 2\kappa\frac{a+\sigma'}{c_0}\Big)^m$$
\BGE \le O_t(1) \exp\Big(\frac2\kappa at+\frac 2\kappa\frac{a+\sigma'}{c_0}(x-3t)-\frac 2\kappa am_0+\frac 2\kappa\frac{a+\sigma'}{c_0}m_0\Big)\le O_t(1)e^{\frac 2\kappa a(x-2t)} \label{Temp2<}\EDE
where the second last inequality follows from (\ref{ac01}), and the last inequality follows from the fact that $|m_0-(3t-x)|<1$.

Suppose ${\cal G}_m\cap{\cal E}_l$ occurs. From (\ref{Em<*}) we have $X_x(s)-2(t+s)\le x-2t+m+(a-2)s$, $0\le s\le S$.
Suppose that $2t\le x\le 3t$. Then $x-2t+m\ge 0$ for any $m\in\N$. Let $p= \frac{x-2t+m}{2-a}\ge 0 $. Then we have
$$\int_0^S \lambda(t+s,X_x(s))ds\le \int_0^p \frac 12 ds+\int_p^\infty 2e^{ x-2t+m+(a-2)s}ds$$\BGE= \frac p2+\int_p^\infty 2e^{(a-2)(s-p)}ds=\frac p2+\frac{2}{2-a}=\frac{x-2t+m+4}{2(2-a)}.\label{est-lambda}\EDE
From (\ref{ac02}), (\ref{PE<*}) and (\ref{El}), we have
$$\sum_{m=1}^\infty\EE\Big[{\bf 1}_{({\cal G}_m\sem {\cal G}_{m-1})\cap{\cal E}_l}\ha\Psi_q(t+S,X_x(S))\exp\Big(\sigma\int_0^S \lambda(t+s,X_x(s))ds\Big)\Big]$$
$$\le O_t(1) \sum_{m=1}^\infty\exp\Big(-\frac 2\kappa(m-1)a+ \sigma  \cdot\frac{x-2t+m+4}{2(2-a)}\Big)$$
\BGE \le O_t(1) e^{\frac{\sigma(x-2t)}{2(2-a)}}\sum_{m=1}^\infty \exp\Big(-\frac 2\kappa a+\frac{\sigma}{2(2-a)}\Big)^m\le O_t(1) e^{\frac 2\kappa a(x-2t)}.\label{Temp3<}\EDE
From (\ref{CE*}), (\ref{Temp2<}), and (\ref{Temp3<}), we have
\BGE \ha\Psi_q(t,x)\le O_t(1) e^{\frac 2\kappa a(x-2t)},\quad 2t\le x\le 3t.\label{2t<x<3t}\EDE

Suppose $0\le x\le 2t$. Let $m_1=\lceil 2t-x\rceil$. Then $x-2t+m\ge 0$ if and only if $m\ge m_1$. If $m\ge m_1$, then (\ref{est-lambda}) still holds. Following the argument of (\ref{Temp3<}),  we get
$$\sum_{m=m_1}^\infty\EE\Big[{\bf 1}_{({\cal G}_m\sem {\cal G}_{m-1})\cap{\cal E}_l}\ha\Psi_q(t+S,X_x(S))\exp\Big(\sigma\int_0^S \lambda(t+s,X_x(s))ds\Big)\Big]$$
$$ \le O_t(1) e^{\frac{\sigma(x-2t)}{2(2-a)}}\sum_{m=m_1}^\infty \exp\Big(-\frac 2\kappa a+\frac{\sigma}{2(2-a)}\Big)^m$$
\BGE \le O_t(1)  e^{\frac{\sigma(x-2t)}{2(2-a)}}\exp\Big(-\frac 2\kappa a+\frac{\sigma}{2(2-a)}\Big)^{m_1}\le O_t(1) e^{\frac 2\kappa a(x-2t)},\label{Temp4<}\EDE
where the last inequality holds because $|m_1-(2t-x)|<1$.

For $m<m_1$, we use the estimation:
$$\int_0^S \lambda(t+s,X_x(s))ds\le \int_0^\infty2e^{x-2t+m+(a-2)s}ds=\frac{2e^{x-2t+m}}{2-a}\quad \mbox{on}\quad {\cal G}_m.$$
From (\ref{El}) we see that, when $m<m_1$, on the event ${\cal G}_m\cap{\cal E}_l$,
$$\ha\Psi_q(t+S,X_x(S))\exp\Big(\sigma\int_0^S \lambda(t+s,X_x(s))ds\Big)$$$$=\exp\Big(O_t(1)e^{(c_0-2)t}+\sigma \frac{e^{x-2t+m}}{2-a}\Big) =1+O_t(1)e^{(c_0-2)t}+O_t(1) e^{x-2t+m},$$
where the last equality holds because $ {e^{x-2t+m}}\le O_t(1)$ for $m<m_1$. Thus,
$$\sum_{m=1}^{m_1-1}\EE\Big[{\bf 1}_{({\cal G}_m\sem {\cal G}_{m-1})\cap{\cal E}_l}\Big(\ha\Psi_q(t+S,X_x(S))\exp\Big(\sigma\int_0^S \lambda(t+s,X_x(s))ds\Big)-1\Big)\Big]$$
$$\le \sum_{m=1}^{m_1-1} e^{-\frac 2\kappa (m-1)a}O_t(1)(e^{(c_0-2)t}+ e^{x-2t+m})=O_t(1)\Big(e^{(c_0-2)t}+ e^{x-2t}\sum_{m=1}^{m_1-1} \exp\Big(1-\frac 2\kappa a\Big)^m\Big)$$
$$\le O_t(1)e^{(c_0-2)t}+O_t(1) e^{x-2t}(1+e^{(1-\frac 2\kappa a)m_1})\le O_t(1)e^{(c_0-2)t}+O_t(1) (e^{x-2t}+e^{\frac 2\kappa a(x-2t)}),$$
where the second last inequality holds because $a\ne \frac \kappa 2$. The above inequality together with (\ref{CE*}), (\ref{Temp2<}), and (\ref{Temp4<}) implies that,
$$\ha\Psi_q(t,x)-1\le O_t(1)(e^{-t}+e^{(c_0-2)t}+e^{x-2t}+e^{\frac 2\kappa a(x-2t)})\le O_t(1) e^{(1-\frac{c_0}2)(x-2t)},\quad 0\le x\le 2t.$$ 
Since  $1\le \ha\Psi_q$, the above inequality implies that
$$ 0\le\ln(\ha\Psi_q(t,x))\le O_t(1) e^{(1-\frac{c_0}2)(x-2t)},\quad 0\le x\le 2t,$$ 
which finishes the proof of (i) for $F=\ha\Psi_q$. The above inequality together with (\ref{2t<x<3t}) implies that
$$0\le \ln(\ha\Psi_q(t,x))\le O_t(1) +\frac 2\kappa a(0\vee (x-2t)),\quad 0\le x\le 3t,$$ 
which finishes the proof of (ii) for $F=\ha\Psi_q$. $\Box$

\section{Annulus SLE with Domain Changed}\label{Sec5}
We now start proving Theorem \ref{main-ann}. The proof will be finished at the end of Section \ref{Section-RN}. Let $p>0$, $\kappa\in(0,4]$, $s\in\R$ $z_0\in\TT$, $w_0\in\TT_p$, and the hull $L$ be as in Theorem \ref{main-ann}. Choose $x_0,y_0\in\R$ such that $z_0=e^{ix_0}$ and $w_0=e^{iy_0-p}$. Let $y_m=y_0+2m\pi$, $m\in\Z$.

Note that $\A_p\sem L$ is a doubly connected domain, whose boundary contain $\TT_p$ and $e^{ix_0}$. Let $p_L=\modd(\A_p\sem L)$. Let $\til L=(e^i)^{-1}(L)$. Then $\til L$ is a subset of $\St_p$ with period $2\pi$.
We may find $W_L$ and $\til W_L$ such that $W_L:(\A_p\sem L;\TT_p)\conf (\A_{p_L};\TT_{p_L})$, $\til W_L:(\St_p\sem \til L;\R_p)\conf(\St_{p_L};\R_{p_L})$, $e^i\circ \til W_L=W_L\circ e^i$, and $\til W_L$ has progressive period $(2\pi;2\pi)$, and

\subsection{Stochastic differential equations}\label{Section-SDE}
Suppose $\xi\in C([0,p))$ with $\xi(0)=x_0$.
Let $g(t,\cdot)$ and $\til g(t,\cdot)$, $0\le t<p$, be the annulus and covering annulus Loewner maps of modulus $p$, respectively, driven by $\xi$. Let $K(t)$ and $\til K(t)$ be the corresponding hulls and covering hulls. Suppose $\xi$ generates a simple annulus Loewner trace $\beta$ of modulus $p$ with $\beta((0,p))\subset \A_p$. Then $\xi$ also generates a simple covering annulus Loewner trace $\til\beta$ of modulus $p$ with $\til\beta((0,p))\subset\St_p$. We have $\beta=e^i\circ\til\beta$, $\til\beta(0)=\xi(0)=x_0$, $K(t)=\beta((0,t])$, and $\til K(t)= \til\beta((0,t])+2\pi\Z$, $0\le t<p$.

Let $T$ be the biggest number in $(0,p]$ such that $\beta((0,T))\cap L=\emptyset$. 
Let $\til\beta_L(t)=\til W_L(\til\beta(t))$ and $\beta_L(t)=W_L(\beta(t))$, $0\le t<T$. Then $\beta_L$ and $\til\beta_L$ are simple curves,  $\beta_L=e^i\circ \til\beta_L$, $\beta_L(0)\in\TT$, and $\beta_L((0,T))\subset\A_{p_L}$. 
Let $v(t)=\ccap_{\A_{p_L}}(\beta_L((0,t]))$.
Then $v$ is a continuous increasing function, which maps $[0,T)$ onto $[0,S)$ for some $S\in(0,p_L]$. Let $\gamma_L(t)=\beta_L(v^{-1}(t))$, $0\le t<S$. Then $\gamma_L(t)$, $0\le t<S$, is the annulus Loewner trace of modulus $p_L$ driven by some $\eta_L\in C([0,S))$. Let $h_L(t,\cdot)$ and $\til h_L(t,\cdot)$, $0\le t<S$, be the  annulus and covering annulus Loewner maps of modulus $p_L$, respectively, driven by $\eta_L$.

For $0\le t<T$, define  $\xi_L(t)=\eta_L(v(t))$, $\til g_L(t,\cdot)=\til h_L(v(t),\cdot)$;
\BGE   \til g_{L,W}(t,\cdot)=\til g_L(t,\cdot)\circ \til W_L;\quad \til W(t,\cdot)=\til g_{L,W}(t,\cdot)\circ g(t,\cdot)^{-1};\label{W(t)}\EDE
Then both $\til g_{L,W}(t,\cdot)$ and $\til W(t,\cdot)$ have progressive period $(2\pi;2\pi)$, and 
\BGE \til g_{L,W}(t,\cdot):(\St_{p}\sem (\til L\cup(\til\beta((0,t])+2\pi\Z));\R_{p})\conf (\St_{p_L-v(t)};\R_{p_L-v(t)}), \label{conf-g-L-til}\EDE
\BGE \til W(t,\cdot):(\St_{p-t}\sem \til L_t;\R_{p-t})\conf(\St_{p_L-v(t)};\R_{p_L-v(t)}),\label{conf-W-til}\EDE
where $\til L_t:=\til g(t,\til L)\subset \St_{p-t}$. We have $g_L(\beta_L(t))=e^{i\xi_L(t)}$. Since $\beta_L=e^i(\til\beta_L)$, there is $n\in\Z$ such that $\til g_L(t, \til\beta_L(t))=\xi_L(t)+2n\pi$ for $0\le t<T$. We now add $2n\pi$ to the driving function $\eta_L$. Then the new $\eta_L$ is still the driving function for $\gamma_L$, $h_L(t,\cdot)$ and $\til h_L(t,\cdot)$, and we have
\BGE \til g_{L,W}(t,\til\beta(t))=\xi_L(t);\label{trace-L}\EDE
 \BGE \til W(t,\xi(t))=\xi_L(t).\label{W(xi)}\EDE

Define $q_m(t)$, $q_{L,m}(t)$, $A_j(t)$, $A_{I,m}(t)$, $X_m(t)$, $X_{L,m}(t)$, $0\le t<T$, such that
\BGE q_m(t)+(p-t)i=\til g(t,y_m+p i);\label{qm}\EDE
\BGE q_{L,m}(t)+(p_L-v(t))i= \til g_{L,W}(t, y_m+p i);\label{qLm}\EDE
\BGE A_j(t)=\til W^{(j)}(t,\xi(t)),\quad j=1,2,3;\quad A_{I,m}(t)=\til W'(t, q_m(t)+(p-t)i);\label{A1I}\EDE
\BGE X_m(t)=\xi(t)-q_m(t);\quad X_{L,m}(t)=\xi_L(t)-q_{L,m}(t).\label{X}\EDE

A standard argument together with Lemma 2.1 in \cite{Zhan} shows that
\BGE v'(t) =\til W'(t,\xi(t))^2=A_1(t)^2.\label{v'*}\EDE
Hence,
\BGE \pa_t \til g_{L,W}(t,z)=\til W'(t,\xi(t))^2 \HA(p_L-v(t), \til g_{L,W}(t,z)-\xi_L(t)).\label{patg-L}\EDE
Since $\HA_I(t,\cdot)$ is odd, from (\ref{ODE-HA-I*}), (\ref{deriv2*}),  and (\ref{patg-L}) we have
\BGE dq_m(t)=-\HA_I(p-t,X_m(t))dt;\label{dq(t)}\EDE
\BGE \frac{d \til g'(t, y_m+p i)}{\til g'(t,y_m+pi)}= \HA_I'(p-t,X_m(t))dt;\label{gI'}\EDE
\BGE  d q_{L,m}(t)=-A_1(t)^2 \HA_I(p_L-v(t),X_{L,m}(t ))dt;\label{dqK}\EDE
\BGE \frac{d \til g_{L,W}'(t, y_m+p i)}{\til g_{L,W}'(t,y_m+pi)}=A_1(t)^2  \HA_I'(p_L-v(t),X_{L,m}(t) )dt.\label{gKI'}\EDE
From (\ref{W(t)}), (\ref{gI'}) and (\ref{gKI'}) we get
 \BGE \frac{dA_{I,m}(t)}{A_{I,m}(t)}=A_1(t)^2 \HA_I'(p_L-v(t),X_{L,m}(t))dt-\HA_I'(p-t,X_m(t))dt.\label{dAI}\EDE

Differentiating $\til W(t,\cdot)\circ \til g(t,z)=\til g_{L,W}(t,z)$ w.r.t.\ $t$ using (\ref{annulus-eqn-covering*}) and (\ref{patg-L}), and letting $w=\til g(t,z)$, we obtain an equality for $\pa_t\til W(t,w)$ with $w\in\St_{p-t}\sem\til L_t$. Differentiating this equality w.r.t.\ $w$, we get an equality for $\pa_t\til W'(t,w)$. Letting $w\to \xi(t)$ in $\St_{p-t}\sem \til L_t$ in these two equalities and using (\ref{Taylor}) we get
\BGE \pa_t\til W(t,\xi(t))=-3A_2(t).\label{-3*}\EDE
 \BGE \frac{\pa_t A_1(t)}{A_1(t)}=\frac 12\Big(\frac{A_2(t)}{A_1(t)}\Big)^2-\frac 43 \frac{A_3(t)}{A_1(t)} +A_1(t)^2\rA(p_L-v(t))- \rA(p-t).\label{patW'}\EDE



Let $\kappa\in (0,4]$. Suppose now $(\xi)$ is a semimartingale, and $d\langle \xi\rangle_t=\kappa dt$, $0\le t<p$. 
We will frequently apply It\^o's formula (c.f.\ \cite{RY}). From (\ref{W(xi)}) and (\ref{-3*})  we have
\BGE d\xi_L(t)=A_1(t)d\xi(t)+\Big(\frac\kappa 2-3\Big) A_2(t)dt.\label{dxiK}\EDE
From (\ref{dq(t)}) and (\ref{dqK})  we see that $X_m(t)$ and $X_{L,m}(t)$ satisfy
\BGE d X_m(t)=d\xi(t)+\HA_I(p-t,X_m(t))dt;\label{dXt}\EDE
\BGE dX_{L,m}(t)=d\xi_L(t)+A_1(t)^2\HA_I(p_L-v(t),X_{L,m}(t))dt.\label{dXK}\EDE
From (\ref{patW'}) we see that
$$\frac{dA_1(t)}{A_1(t)}=\frac{A_2(t)}{A_1(t)}d\xi(t)+\Big[\frac 12\Big(\frac{A_2(t)}{A_1(t)}\Big)^2+\Big(\frac \kappa 2-\frac 43\Big) \frac{A_3(t)}{A_1(t)} +A_1(t)^2\rA(p_L-v(t))- \rA(p-t)\Big]dt.$$
Let $\cc$ and $\alpha$ be as in (\ref{cc}) and (\ref{alpha}), respectively. Then we compute that
\BGE \frac{dA_1(t)^\alpha}{A_1(t)^\alpha}=\alpha \frac{A_2(t)}{A_1(t)}d\xi(t)+\Big[\frac{\cc}6A_S(t)+\alpha A_1(t)^2\rA(p_L-v(t)) - \alpha\rA(p-t)\Big] dt,\label{dA1}\EDE
where $A_S(t):= \frac{A_3(t)}{A_1(t)}-\frac 32(\frac{A_2(t)}{A_1(t)})^2$  is the Schwarz derivative of $\til W(t,\cdot)$ at $\xi(t)$.

Let
\BGE Y_m(t)=\Gamma_0(p-t,X_m(t)),\quad Y_{L,m}(t)=\Gamma_0(p_L-v(t),X_{L,m}(t)),\label{Y=X-m}\EDE
From (\ref{Lambda-Gamma}),  (\ref{PDE-Gamma}),  (\ref{v'*}), (\ref{dXt}) and (\ref{dXK}) we find that
\BGE \frac{dY_m(t)}{Y_m(t)}=\frac 1\kappa {\Lambda_0(p-t,X_m(t))}d\xi(t)-\alpha \HA_I'(p-t,X_m(t))dt,\label{dY-m}\EDE
\BGE \frac{dY_{L,m}(t)}{Y_{L,m}(t)}=\frac 1\kappa {\Lambda_0(p_L-v(t),X_{L,m}(t))}d\xi_L(t)-\alpha A_1(t)^2 \HA_I'(p_L-v(t),X_{L,m}(t))dt.\label{dYK-m}\EDE

Define $M_m$ on $[0,T)$ by
\BGE M_m=A_1^\alpha A_{I,m}^\alpha \frac{Y_{L,m}}{Y_m}\exp\Big(-\frac{\cc}6\int_0^\cdot A_S(s)ds+\alpha\int_{p-\cdot}^{p_L-v(\cdot)} \rA(s)ds\Big).\label{Mm}\EDE
From (\ref{dAI}), (\ref{dxiK}), (\ref{dA1}), (\ref{dY-m}) and (\ref{dYK-m}), we find that
$$\frac{dM_m(t)}{M_m(t)}=\Big[\alpha \frac{A_2(t)}{A_1(t)} +\frac {A_1(t)}\kappa {\Lambda_0(p_L-v(t),X_{L,m}(t))} -\frac 1\kappa {\Lambda_0(p-t,X_m(t))}\Big]\cdot $$
\BGE \cdot(d\xi(t)-\Lambda_0(p-t,X_m(t))dt),\quad 0\le t<T.\label{dM/M}\EDE

Let $C_{p,L}=\exp(-\frac{\cc}2\int_p^{p_L}(\rA(s)+\frac 1s)ds)>0$. We have $M_m=N_m\exp(\cc U)$, where 
\BGE N_m= C_{p,L}{A_1^\alpha A_{I,m}^\alpha} \frac{Y_{L,m}}{Y_m} \exp\Big((\alpha+\frac{\cc}2)\int^{p_L-v(\cdot)}_{p-\cdot}\rA(s)ds+\frac{\cc}2\int^{p_L-v(\cdot)}_{p-\cdot}\frac 1sds\Big);\label{N}\EDE
\BGE U=-\frac 16\int_0^\cdot A_S(s)ds-\frac12\int^{p_L-v(\cdot)}_{p-\cdot}(\rA(s)+\frac 1s)ds+\frac 12\int^{p_L}_p (\rA(s)+\frac 1s)ds.\label{U(t)}\EDE

\subsection{Rescaling} \label{Section-Rescaling}
Let $\ha p=\frac{\pi^2}p$, $\ha T=\frac{\pi^2}{p-T}-\ha p$ and $\check t=p-\frac{\pi^2}{\ha p+t}$.
Then the function $t\mapsto \check t$ maps $[0,\infty)$ onto $[0,p)$, and maps $[0,\ha T)$ onto $[0,T)$.
Let $\ha L=\frac \pi p\til L\subset\St_\pi$, $\ha x_0=\frac \pi p x_0$, and $\ha y_m=\frac\pi p y_m$, $m\in\Z$. Since $\til L$ has period $2\pi$, and $\dist(\til L,\{x_0\}\cup \R_p)>0$, we see that $\ha L$ has period $2\ha p$, and $\dist(\ha L,\{\ha x_0\}\cup \R_\pi)>0$. Let $\ha \beta(t)=\frac\pi p \til \beta (\check t)$, $0\le t<\infty$. Then $\ha\beta$ is a curve in $\St_\pi\cup\R$ started from $\ha x_0$, and $\ha T$ is the biggest number in $(0,\infty]$ such that $\ha\beta((0,\ha T))\cap \ha L=\emptyset$. Furthermore, we have
\BGE \{T=p\}=\{\ha T=\infty\}=\{\beta\cap L=\emptyset\}=\{\til\beta\cap \til L=\emptyset\}=\{\ha\beta\cap\ha L=\emptyset\}.\label{T=p}\EDE
 For $0\le t<\ha T$, let $\ha p_{L}(t)=\frac{\pi^2}{p_L-v(\check t)}$. Since $p_L-v(t)=\modd(\A_p\sem L\sem \beta((0,t])$, while $p-t=\modd(\A_p\sem \beta((0,t])$, we have $p_L-v(t)\le p-t$, $0\le t<T$. Thus,
\BGE \ha p_{L}(t)\ge \frac{\pi^2}{p-\check t}= \ha p+t,\quad 0\le t<\ha T.\label{pL>p}\EDE
 From (\ref{rA-ha}), for any $0\le t<\ha T$,
\BGE -\int_{\ha p+t}^{\ha p_{L}(t)} \ha\rA(s)ds=\int_{p-\check t}^{{p_{ L}-v(\check t)}} \Big(\rA(  s)+\frac 1{  s}\Big)d s .\label{int-rA-ha=}\EDE

Let $m\in\Z$. For $0\le t<\ha T$, define
\BGE\ha\xi(t)=\frac{\ha p+t}{\pi}\cdot \xi (\check t);\quad \ha q_m(t)=\frac{\ha p+t}{\pi}\cdot q_m(\check t);\quad \ha X_m(t)= \frac{\ha p+t}{\pi}\cdot X_m(\check t);\label{haX}\EDE
\BGE\ha\xi_{L}(t)=\frac{\ha p_{L}(t)}\pi\cdot\xi_L(\check t);\quad \ha q_{L,m}(t)=\frac{\ha p_{L}(t)}\pi\cdot q_{L,m}(\check t); \quad \ha X_{L,m}(t)=\frac{\ha p_{L}(t)}\pi\cdot X_{L,m}(\check t).\label{haXL}\EDE
From (\ref{X}) we have $\ha X_m=\ha\xi-\ha q_m$ and $\ha X_{L,m}=\ha\xi_L-\ha q_{L,m}$. For $0\le t<\ha T$, define
\BGE\ha g(t,z)=\frac{\ha p+t}\pi\,\til g({\check t},\frac p\pi z);\quad \ha g_{L,W}(t,z)=\frac{\ha p_L(t)}{\pi}\,\til g_{L,W}(\check t,\frac p\pi z).\label{ha-g-L}\EDE
From (\ref{conf-g-til}), (\ref{trace-til*}), (\ref{trace-L}), (\ref{conf-g-L-til}), (\ref{qm}) and (\ref{qLm}) we have 
\BGE \ha g(t,\cdot):(\St_\pi\sem(\ha \beta((0,t])+2\ha p\Z); \R_\pi,\ha\beta(t),\ha y_m+\pi i)\conf(\St_{\pi};\R_\pi,\ha\xi(t),\ha q_m(t) );\label{conf-g-ha}\EDE
\BGE \ha g_{L,W}(t,\cdot):(\St_\pi\sem ((\ha \beta((0,t])+2\ha p\Z)\cup\ha L); \R_\pi,\ha\beta(t),\ha y_m+\pi i)\conf(\St_{\pi};\R_\pi ,\ha\xi_L(t),\ha q_{L,m}(t)+\pi i).\label{conf-g-L-ha}\EDE
For $0\le t<\ha T$, define
\BGE \ha W(t,\cdot)=\ha g_{L,W}(t,\cdot)\circ \ha g(t,\cdot)^{-1};
\label{W-ha}\EDE
\BGE \ha A_1(t)=\ha W'(t,\ha\xi(t)),\quad \ha A_{I,m}(t)=\ha W'(t, \ha q_m(t)+\pi i).\label{ha-A-def}\EDE
Let $\ha L_t=\frac{\ha p+t}{\pi} \til L_{\check t}$. Since $\til L_t=\til g(t,\til L)$, we have $\ha L_t=\ha g(t,\ha L)$.
From  (\ref{conf-g-ha}) and (\ref{conf-g-L-ha}) we have
\BGE \ha W(t,\cdot) :(\St_\pi\sem \ha L_t;\R_\pi)\conf (\St_\pi;\R_\pi);\label{haW=}\EDE
\BGE \ha W(t,\ha\xi(t))=\ha\xi_{L}(t);\quad \ha W(t,\ha q_m(t)+\pi i)=\ha q_{L,m}(t)+\pi i.\label{W(q)-ha-m}\EDE
From  (\ref{W(t)}), (\ref{A1I}), (\ref{ha-g-L}) and (\ref{W-ha}) we have
 \BGE \ha A_1(t)=\frac{\ha p_{L}(t)}{\ha p+t}\,A_1(\check t),\quad   \ha A_{I,m}(t)=\frac{\ha p_{L}(t)}{\ha p+t}\, A_{I,m}(\check t).\label{A-ha=}\EDE

 For $m\in\Z$ and $0\le t<\ha T$, let
\BGE \ha Y_m(t)=\ha\Gamma_0(\ha p+t,\ha X_m(t)),\quad \ha Y_{{L},m}(t)=\ha\Gamma_0(\ha p_{L}(t),\ha X_{L}(t)).\label{Y=X-m-ha}\EDE
From (\ref{ha-Gamma}),  (\ref{Y=X-m}), (\ref{haX}), and (\ref{haXL}),   we have
 \BGE \ha Y_m(t)=\Big(\frac\pi{\ha p+t}\Big)^\alpha Y_m(\check t),\quad \ha Y_{{L},m}(t)=\Big(\frac\pi{\ha p_{L}(t)}\Big)^\alpha Y_{L,m}(\check t).\label{Y-ha=}\EDE

Define $\ha N_m$ on $[0,\ha T)$ such that
$$\ha N_m=C_{p,L}\ha A_1^\alpha \ha A_{I,m}^\alpha {\ha Y_{{L},m}}{\ha Y_m}^{-1}\exp\Big(-(\alpha+\frac{\cc}2)\int_{\ha p+\cdot}^{\ha p_{L}(\cdot)} \ha \rA( s)d s\Big).$$
 From (\ref{N}),  (\ref{int-rA-ha=}),   (\ref{A-ha=}) and  (\ref{Y-ha=}) we find that
 \BGE \ha N_m(t)=N_m(\check t),\quad 0\le t<\ha T.\label{N-N}\EDE
From  (\ref{cc}), (\ref{RA-rA-ha}),  (\ref{sigma-tau}), (\ref{alpha}), (\ref{Gamma-m-ha}), (\ref{ha-Gamma-infty}), and (\ref{Y=X-m-ha}), we see that for $0\le t<\ha T$,
$$ \ha N_m(t)=\ha A_1( t)^\alpha \ha A_{I,m}( t)^\alpha  \, \ha \Gamma_q(\ha p_{L}(t), \ha X_{{L},m}(t))\,\ha\Gamma_q(\ha p+t,\ha X_m(t))^{-1}\cdot$$ \BGE\cdot\exp\Big(-\alpha\int_{\ha X_m( t)}^{\ha X_{{L},m}( t)} \tanh_2(s)ds-(\alpha+\frac{\cc}2)(\ha\RA(\ha p_{L}(t))-\ha\RA(\ha p+t))\Big).\label{N-bd}\EDE

\section{Estimations on $\ha N_m(t)$} \label{section-est}
For $m\in\Z$, 
let ${\cal P}_m$ denote the set of $(\rho_1,\rho_2)$ with the following properties.
\begin{enumerate}
  \item For $j=1,2$, $\rho_j$ is a polygonal crosscut in $\St_\pi$ that grows from a point on $\R$ to a point on $\R_\pi$, whose line segments are  parallel to either $x$-axis or $y$-axis, and whose vertices other than the end points  have rational coordinates.
  \item $\rho_1+2n\ha p$, $\rho_2+2n\ha p$, $n\in\Z$, and $\ha L$ are mutually disjoint; $\rho_1$ lies to the left of $\rho_2$.
  \item $\rho_1\cup\rho_2$ disconnects $\ha x_0$ and $\ha y_m+ \pi i$ from $\ha L$ in $\St_\pi$.
\end{enumerate}

For each $( \rho_1, \rho_2)\in{\cal P}_m$, let $\ha T_{\rho_1,\rho_2}$  denote the biggest number such that $\ha\beta((0,\ha T_{\rho_1,\rho_2}))\cap(\rho_1\cup\rho_2)=\emptyset$. Since $\ha\beta$ starts from $\ha x_0$, we have $\ha T_{\rho_1,\rho_2}\le \ha T$. Let ${\cal E}_m$ be as in Section \ref{section-limit}. Then
$${\cal E}_m=\{\lim_{t\to p}\til\beta(t)=y_m+p i\}=\{\lim_{t\to\infty}\ha\beta(t)=\ha y_m+\pi i\}.$$
We will prove the following proposition  at the end of Section \ref{conformal-section} and  Section \ref{compact-section}.

\begin{Proposition} Let $m\in\Z$.
\begin{enumerate}
    \item [(i)]  $\lim_{t\to \infty} \ln( \ha N_m(t)/C_{p,L})=0)$ on the event ${\cal E}_m\cap\{\ha T=\infty\}$.
  \item [(ii)] For any $(\rho_1,\rho_2)\in {\cal P}_m$, $\ln(\ha N_m(t))$ is uniformly bounded on $[0,\ha T_{\rho_1,\rho_2})$.
\end{enumerate}
 \label{N-to-1-ha}
\end{Proposition}

For $m\in\Z$, define $\til{\cal P}_m=\{(\frac p\pi\rho_1,\frac p\pi\rho_2):(\rho_1,\rho_2)\in{\cal P}_m\}$.
Then for each $(\rho_1,\rho_2)\in\til{\cal P}_m$, $\rho_1$ and $\rho_2$ are simple  curves that grow from a point on $\R$ to a point on $\R_\pi$, and $\rho_1\cup\rho_2$ disconnects $x_0$ and $y_m+p i$ from $\til L$. For each $(\rho_1,\rho_2)\in\til{\cal P}_m$, let $T_{\rho_1,\rho_2}$ denote the biggest time such that $\til\beta((0,T_{\rho_1,\rho_2}))\cap(\rho_1\cup\rho_2)=\emptyset$. Then $T_{\rho_1,\rho_2}\le T$, and the function $t\mapsto \check t$ maps $[0,\ha T_{\rho_1,\rho_2})$ onto $[0,T_{\frac p\pi\rho_1,\frac p\pi\rho_2})$.
From Proposition \ref{N-to-1-ha},  (\ref{T=p}), and (\ref{N-N}) we conclude the following proposition:

\begin{Proposition}
  Let $m\in\Z$.
\begin{enumerate}
    \item [(i)]  $\lim_{t\to p}  N_m(t)=C_{p,L}$ on the event ${\cal E}_m\cap\{T=p\}$.
  \item [(ii)] For any $(\rho_1,\rho_2)\in \til{\cal P}_m$, $N_m(t)$ is uniformly bounded on $[0,T_{\rho_1,\rho_2})$.
\end{enumerate}
 \label{N-to-1}
\end{Proposition}

\subsection{The limit value} \label{conformal-section}
We now use $\HH$ and $\D$ to denote the upper half-plane $\{\Imm z>0\}$ and the unit disk $\{|z|<1\}$, respectively. Let ${\cal H}$ denote the set of bounded hulls in $\HH$. For every $H\in\cal H$, there is unique $\vphi_H$ which maps $\HH\sem H$ conformally onto $\HH$ such that as $z\to\infty$, $\vphi_H(z)=z+\frac cz+O(|z|^{-2})$, where $c=:\hcap(H)$ is called the half-plane capacity of $H$. If $H=\emptyset$, then $\vphi_H=\id$ and $\hcap(H)=0$; otherwise $\hcap(H)>0$.

Suppose $H\in\cal H$ and $H\ne\emptyset$. Then
$\lin{H}\cap\R\ne\emptyset$. Let $a_H=\inf(\lin{H}\cap\R)$ and
$b_H=\sup(\lin{H}\cap\R).$ Let
$$\Sigma_H=\C\sem(H\cup\{\lin{z}:z\in
H\}\cup[a_H,b_H]).$$ By reflection principle, $\vphi_H$ extends to
$\Sigma_H$, and $\vphi_H: \Sigma_H\conf\C\sem[c_H,d_H]$
for some $c_H<d_H\in\R$. Moreover, $\vphi_H$ is increasing on
$(-\infty,a_H)$ and $(b_H,\infty)$, and maps them onto
$(-\infty,c_H)$ and $(d_H,\infty)$, respectively. So $\vphi_H^{-1}$
extends conformally to $\C\sem[c_H,d_H]$.
\vskip 4mm

\begin{Example} Suppose $r>0$. Let $H=\{z\in\HH:|z|\le r\}$. Then $H\in\cal H$, $a_H= -r$ and $b_H=r$. It is clear that $\vphi_H(z)=z+\frac{r^2}{z}$. Thus $\hcap(H)=r^2$, and $[c_H,d_H]=[-2r,2r]$. \label{Example}\end{Example}

From (5.1) in \cite{LERW} there is a measure $\mu_H$ supported by $[c_H,d_H]$ with $|\mu_H|=\hcap(H)$ such that for any $z\in\Sigma_H$,
\begin{equation}\vphi_{H}^{-1}(z)-z=\int_{c_H}^{d_H}\frac{-1}{z-x}\,
d\mu_H(x). \label{chordal equation general} \end{equation}
Since $\vphi_\emptyset=\id$, (\ref{chordal equation general}) is also true for $H=\emptyset$ if we set $\mu_\emptyset=0$, $a_\emptyset=c_\emptyset=\infty$ and $b_\emptyset=d_\emptyset=-\infty$. The following lemma is a combination of Lemma 5.2 and Lemma 5.3 in \cite{LERW}.
\begin{Lemma} \begin{enumerate} \item [(i)] For any $H\in\cal H$, $\vphi_H(x)\le x$ on $(-\infty,a_H)$ and $\vphi_H(x)\ge x$ on $(b_H,\infty)$;
 \item [(ii)] If $H_1,H_2\in\cal H$ and $ H_1\subset H_2$, then $[c_{H_1},d_{H_1}]\subset [c_{H_2},d_{H_2}]$.
  \end{enumerate}
  \label{comb-Lem}
\end{Lemma}

\begin{Lemma}
  Let $r>0$. Suppose $H\in\cal H$ and $H\subset\{|z|\le r\}$. Suppose $W:(\HH\sem H;\infty)\conf(\HH;\infty)$,  and satisfies that $W'(\infty)=1$, $W((-\infty,-r))\subset(-\infty,0)$ and $W((r,\infty))\subset(0,\infty)$. Then for any $z\in\HH\cup\R$ with $|z|\ge (1+r)^2$, $|W(z)-z|\le r^2+2r$. \label{1}
\end{Lemma}
{\bf Proof.} Let $H_r=\{z\in\HH:|z|\le r\}$. Then $H\subset H_r\in\cal H$. From Lemma \ref{comb-Lem} (i) and Example \ref{Example}, we have $[c_H,d_H]\subset [c_{H_r},d_{H_r}]=[-2r,2r]$. Define $K=\vphi_H(H_r\sem H)$. Then $K\in\cal H$ and $\vphi_{H_r}=\vphi_K\circ \vphi_H$, which implies that $\hcap(H_r)=\hcap(K)+\hcap(H)$. Thus, $\hcap(H)\le \hcap(H_r)=r^2$. Applying Lemma \ref{comb-Lem} (ii) to $K$, we find that $\vphi_{H}(x)\le \vphi_{H_r}(x)$ for any $x\in (b_{H_r},\infty)=(r,\infty)$. Thus, $\inf \vphi_H((r,\infty))\le \inf \vphi_{H_r}((r,\infty))=d_{H_r}=2r$. Similarly, $\sup\vphi_H((-\infty,-r))\ge -2r$. Since both $W$ and $\vphi_H$ map $\HH\sem H$ conformally onto $\HH$ and fix $\infty$, and have derivative $1$ at $\infty$, there is $w\in\R$ such that $W(z)=\vphi_H(z)-w$ for any $z\in\HH\sem H$. From the assumption of $W$, we get $\inf W((r,\infty))\ge 0\ge \sup W((-\infty,-r))$. So we have $|w|\le 2r$. We now suffice to show that for any $z\in\HH\cup\R$ with  $|z|\ge (1+r)^2$, $|\vphi_H(z)-z|\le r^2$.

Let $z\in\HH\cup\R$ and $|z|\ge 1+2r$. Since $[c_H,d_H]\subset [-2r,2r]$ and $|\mu_H|=\hcap(H)\le r^2$, from (\ref{chordal equation general}) we have $|\vphi_H^{-1}(z)-z|\le r^2$. Let $\gamma=\{z\in\HH:|z|=1+2r\}$. Then $\vphi_H^{-1}(\gamma)$ is a crosscut in $\HH$, which divides $\HH$ into two components. Let $D$ denote the unbounded component. Then $\vphi_H^{-1}$ maps $\{z\in\HH\cup\R:|z|\ge 1+2r\}$  onto $\lin{D}$.
Since $|\vphi_H^{-1}(z)-z|\le r^2$ for $z\in\gamma$, we have $\vphi_H^{-1}(\gamma)\subset\{|z|\le 1+2r+r^2\}$, which implies that $\lin{D}\supset\{z\in\HH\cup\R:|z|\ge (1+r)^2\}$. Thus, if $z\in\HH\cup\R$ and $|z|\ge (1+r)^2$, then $\vphi_H(z)\in \{z\in\HH\cup\R:|z|\ge 1+2r\}$, which implies that $|\vphi_H(z)-z|=|\vphi_H^{-1}(\vphi_H(z))-\vphi_H(z)|\le r^2$. $\Box$

\begin{Lemma}
  Let $K$ be a hull in $\St_\pi$  such that $\Ree z\le c$ for any $z\in K$. Suppose that $V:(\St_\pi\sem K;+\infty)\conf(\St_\pi;+\infty)$, and satisfies $V((c, \infty))\subset\R$ and $V((c,  \infty)+\pi i)\subset \R_\pi$. Then there is $h\in\R$ such that if $z\in\St_\pi\cup\R\cup\R_\pi$ and $\Ree z\ge c+\ln(4)$, then $|V(z)-z-h|\le 12e^{c-\Ree z}$. \label{2}
\end{Lemma}
{\bf Proof.} There is $h\in\R$ such that $V(z)=z+h+o(1)$ as $z\in\St_\pi$ and $z\to+\infty$. By considering $V-h$ instead of $V$, we may assume that $V(z)=z+o(1)$ as $z\in\St_\pi$ and $z\to+\infty$. Let $z\in\St_\pi\cup\R\cup\R_\pi$ with $\Ree z\ge c+\ln(4)$. Let $a=\Ree z-c$. Then $e^a\ge 4$, and there is $r\in(0,1]$ such that $(1+r)^2/r=e^a$. Let $H=\frac{r}{e^c}\exp(K)$ and \BGE W(z)=\exp(V(\ln z-\ln(r)+c)+\ln(r)-c).\label{W}\EDE Then $H\in\cal H$, $H\subset\{z\in\HH:|z|\le r\}$,  $W:(\HH\sem H;\infty)\conf (\HH;\infty)$, and $W'(\infty)=1$. Since $V((c,\infty))\subset\R$ and $V((c,\infty)+\pi i)\subset \R_\pi$, we have $W((-\infty,-r))\subset(-\infty,0)$ and $W((r,\infty))\subset(0,\infty)$. Since $z\in\St_\pi\cup\R\cup\R_\pi$ and $\Ree z=c+a=c+2\ln(1+r)-\ln(r)$, we have $e^{z+\ln(r)-c}\in\HH\cup\R$ and $|e^{z+\ln(r)-c}|\ge (1+r)^2$. From  lemma \ref{1}, we have  $|W(e^{z+\ln(r)-c})-e^{z+\ln(r)-c}|<r^2+2r$. So the line segment $[e^{z+\ln(r)-c},W(e^{z+\ln(r)-c})]$ lies outside $\D$. Since $|\ln'(z)|\le 1$ for $z\in\C\sem\D$, we have $|\ln(W(e^{z+\ln(r)-c}))-({z+\ln(r)-c})|<r^2+2r\le 3r$. From (\ref{W}), $V(z)=\ln(W(e^{z+\ln(r)-c}))-\ln(r)+c$, so $|V(z)-z|<3r$. Finally, since $e^ar=(r+1)^2=r^2+2r+1\le 3r+1$, we have $r\le \frac{1}{e^a-3}$. Since $a\ge \ln(4)$, $e^a-3\ge1\ge 4e^{-a}$.
Thus, $|V(z)-z|\le 12e^{-a}=12e^{c-\Ree z}$. $\Box$

\begin{Proposition}
  Let $K=K_+\cup K_-\subset\St_\pi$, where $\Ree K_-$ is bounded above by $c_-\in\R$, $\Ree K_+$ is bounded below by $c_+\in\R$, and $c_+- c_-\ge 2\ln(12)$. Suppose $W$ maps $\St_\pi\sem K$ conformally onto $\St_\pi$, and satisfies $W((c_-,c_+))\subset\R$ and $W((c_-,c_+ )+\pi i)\subset \R_\pi$. Then there exists $h\in\R$ such that if $z\in\St_\pi$ satisfies that $d:=\min\{c_+-\Ree z,\Ree z-c_-\}\ge \ln(12)$, then $|W(z)-z-h|\le 48 e^{-d}$. \label{ThmK}
\end{Proposition}
{\bf Proof.} Choose $W_-:(\St_\pi\sem K_-;+\infty)\conf (\St_\pi;+\infty)$. By composing a suitable $U:(\St_\pi;+\infty)\conf(\St_\pi;+\infty)$ on its left, we may assume that $W_-$ satisfies $W_-\circ W^{-1}(-\infty)=-\infty$. Let $K_+'=W_-(K_+)$ and $W_+=W\circ W_-^{-1}$. Then $W_+:(\St_\pi\sem K_+';-\infty)\conf(\St_\pi;-\infty)$.

The assumption on $W$ implies that there is no $x\in(c_-,\infty)$ such that $W(x)=-\infty$. Since $W_-\circ W^{-1}(-\infty)=-\infty$, there is no $x\in (c_-,\infty)$ such that $W_-(x)=-\infty$. This implies that $W_-((c_-,\infty))\subset\R$. Similarly, $W_-((c_-,\infty)+\pi i)\subset\R_\pi$. From Lemma \ref{2}, there is $h_-\in\R$ such that
\BGE |W_-(z)-z-h_-|\le 12e^{c_--\Ree z},\quad \mbox{if }z\in\St_\pi\cup\R\cup\R_\pi\mbox{ and } \Ree z\ge c_-+\ln(4).\label{W-z-}\EDE

Let $c_d=c_+-c_-\ge 2\ln(12)$ and $c_+'=c_++h_--12e^{-c_d}$. If $z\in\St_\pi\cup\R\cup\R_\pi$ and $\Ree z\ge c_+$, then $\Ree z\ge c_-+c_d\ge c_-+\ln(4)$, which implies that $\Ree W_-(z)\ge  c_+'$ by (\ref{W-z-}). Since $\Ree K_+$ is bounded below by $c_+$, we find that $\Ree K_+'$ is bounded below by $c_+'$. Suppose $W_-((c_-,c_+))=(a_-,a_+)$. The above argument show that $a_+ \ge c_+'$. Since $W=W_+\circ W_-$ and $W((c_-,c_+))\subset\R$, we have $W_+((a_-,a_+))\subset\R$. Since $W_+$ fixes $-\infty$, we have $W_+((-\infty,a_+))\subset\R$, which implies that $W_+((-\infty,c_+))\subset\R$ as $a_+\ge c_+'$. Similarly, $W_+((-\infty,c_+)+\pi i)\subset\R_\pi$. From a mirror result of Lemma \ref{2}, we see that there exists $h_+\in\R$ such that
\BGE |W_+(w)-w-h_+|\le 12e^{\Ree w-c_+'},\quad\mbox{if }w\in\St_\pi\cup\R\cup\R_\pi\mbox{ and } \Ree w\le c_+'-\ln(4).\label{W-z+}\EDE

Now suppose that $z\in\St_\pi$ satisfies that $d:=\min\{c_+-\Ree z,\Ree z-c_-\}\ge \ln(12)$. From (\ref{W-z-}) and that $c_+-\Ree z,\Ree z-c_-\ge d\ge \ln(12)$, we obtain
$$\Ree W_-(z)\le \Ree z+h_-+12e^{c_--\Ree z}\le c_+-d+h_-+12e^{-d}=c_+'-d+12e^{-c_d}+12e^{-d}$$
\BGE\le c_+'-d+12e^{-2\ln(12)}+12e^{-\ln(12)}=c_+'-d+\frac{13}{12}\le c_+'-d+\ln(3)\le c_+'-\ln(4).\label{d-ln}\EDE
Let $h=h_++h_-$. Applying (\ref{W-z+}) to $w=W_-(z)$ and using (\ref{W-z-}) and (\ref{d-ln}), we get
$$|W(z)-z-h|\le |W_+(W_-(z))-W_-(z)-h_+|+|W_-(z)-z-h_-|$$
$$\le 12 e^{\Ree W_-(z)-c_+'}+12e^{c_--\Ree z}<12e^{\ln(3)-d}+12e^{-d}=48e^{-d}.\quad\Box$$

\vskip 4mm
Differentiating (\ref{chordal equation general}) w.r.t.\ $z$, we see that for $z\in\Sigma_H$,
$$(\vphi_{H}^{-1})'(z)-1=\int_{c_H}^{d_H}\frac{1}{(z-x)^2}\,
d\mu_H(x);\quad
(\vphi_{H}^{-1})^{(n)}(z)=\int_{c_H}^{d_H}\frac{(-1)^{n+1}n!}{(z-x)^{n+1}}\,
d\mu_H(x),\quad n\ge 2.$$
The proofs of Lemma \ref{1}, Lemma \ref{2}, and Proposition \ref{ThmK} can be slightly modified to prove the following proposition.

\begin{Proposition}
  There are constants $C_1,C_2>0$ such that the following hold. Let $K=K_+\cup K_-\subset\St_\pi$, where $\Ree K_+$ is bounded above by $c_-\in\R$, $\Ree K_+$ is bounded below by $c_+\in\R$, and $c_+\ge c_-+2C_2$. Suppose $W$ maps $\St_\pi\sem K$ conformally onto $\St_\pi$, and satisfies $W((c_-,c_+))\subset\R$ and $W((c_-,c_+)+\pi i)\subset \R_\pi$. Then for any $z\in\St_\pi$ with $d:=\min\{c_+-\Ree z,\Ree z-c_-\}\ge C_2$, we have $|W'(z)-1|,|W''(z)|,|W'''(z)|\le C_1 e^{-d}$. \label{ThmK'}
\end{Proposition}

\no {\bf Proof of Proposition \ref{N-to-1-ha} (i).} From (\ref{N-bd}), we suffice to show that (i) holds if $\ln(\ha N_m(t)/C_{p,L})$ is replaced by $\ln(\ha A_1(t))$, $\ln(\ha A_{I,m}(t))$, $\ln(\ha\Gamma_q(\ha p+t,\ha X_m(t)))$, $\ln(\ha\Gamma_q(\ha p_{L}(t),\ha X_{{ L},m}(t)))$, $\ha X_{{L},m}(t)-\ha X_m(t)$, or $\ha \RA(\ha p_{L}(t))-\ha \RA(\ha p+t)$, respectively. Suppose ${\cal E}_m\cap\{\ha T=\infty\}$ occurs. From (\ref{pL>p}) and (\ref{RA-rA-ha}) we conclude that $\ha \RA(\ha p_{L}(t))-\ha \RA(\ha p+t)\to 0$ as $t\to\infty$.

Decompose $\ha L$ into $\ha L_l$ and $\ha L_r$ such that $\lin{\ha L_l\cap\R}$ (resp.\ $\lin{\ha L_r\cap\R}$) lies to the left (resp.\ right) of $\ha x_0$. Let $\ha L_{l,t}=\ha g(t,\ha L_l)$ and $\ha L_{r,t}=\ha g(t,\ha L_r)$, $0\le t<T$. Then $\lin{\ha L_{l,t}\cap\R}$ (resp.\ $\lin{\ha L_{r,t}\cap\R}$) lies to the left (resp.\ right) of $\ha \xi(t)$.
Let $E_r$ be a subset of $\St_\pi$, which touches both $\R$ and $\R_\pi$, disconnects $\ha \beta$ from $\ha L_r$ in $\St_\pi$, and is disjoint from $\ha L$ and $\ha \beta$. As $t\to\infty$, the diameter of $\ha \beta((t,\infty))$ tends to $0$, which implies that the extremal distance (c.f.\ \cite{Ahl}) in $\St_\pi\sem (\ha \beta((0,t])+2\ha p\Z)$ between $E_r$ and the set
$$S_t:=(-\infty,x_0]\cup (\ha \beta((0,t])-2\ha p\N)\cup\{\mbox{the left side of }\ha \beta((0,t])\}\cup \{y+\pi i: y\le \ha y_m\}$$
tends to $\infty$. From (\ref{conf-g-ha}) and conformal invariance, the extremal distance in $\St_\pi $ between $\ha g(t,E_r)$ and $(-\infty,\ha\xi(t)]\cup\{x+\pi i: x\le \ha q_m(t)\}$ tends to $\infty$ as $t\to \infty$. Since $E_r$ touches both $\R$ and $\R_\pi$, $\ha g(t,E_r)$ also has this property. Thus, $\dist(\{\ha\xi(t),\ha q_m(t)\},\ha g(t,E_r))\to \infty$ as $t\to \infty$. Since $E_r$ disconnects $\ha \beta$ from $\ha L_r$ in $\St_\pi$, we see that $\ha g(t,E_r)$ disconnects $\ha\xi(t)$ and $\ha q_m(t)+\pi i$ from $\ha L_{r,t}$. Thus, $\dist(\{\ha\xi(t),\ha q_m(t)+\pi i\},\ha L_{r,t})\to \infty$ as $t\to \infty$. Similarly, $\dist(\{\ha\xi(t),\ha q_m(t)+\pi i\},\ha L_{l,t})\to \infty$ as $t\to \infty$. Thus,  $\dist(\{\ha\xi(t),\ha q_m(t)+\pi i\},\ha L_t\}\to \infty$ as $t\to\infty$. From (\ref{ha-A-def}), (\ref{haW=}), and Proposition \ref{ThmK'} we conclude that $\ln(\ha A_1(t))\to 0$ and $\ln(\ha A_{I,m}(t))\to 0$ as $t\to\infty$.  From (\ref{W(q)-ha-m}) and Proposition \ref{ThmK}, we see that $\ha X_{{L},m}(t)-\ha X_m(t)=(\ha \xi_{L}(t)-\ha \xi(t))-(\ha q_{{ L},m}(t)-\ha q_m(t))\to 0$ as $t\to \infty$.

Let $a_r(t)=\min\{\lin{\ha L_{r,t}}\cap\R\}$ and $a_l(t)=\max\{\lin{\ha L_{l,t}}\cap\R\}$. From the last paragraph, we see that $a_r(t)-\ha\xi(t)$, $a_r(t)-\ha q_m(t)$, $\ha\xi(t)-a_l(t)$, and $\ha q_m(t)-a_l(t)$ all tend to $+\infty$ as $t\to \infty$. Since $\ha X_m=\ha\xi-\ha q_m$, we have $a_r(t)-a_l(t)\pm X_m(t)\to \infty$ as $t\to \infty$. Since $\ha L_t$ has period $2(\ha p+t)$, we have $a_r(t)- a_l(t)\le 2(\ha p+t)$.
Thus, $2(\ha p+t)-|X_m(t)|\to \infty$ as $t\to \infty$. Let $b_1,b_2\in\R$ be such that $b_1<\ha y_m<b_2=b_1+2\ha p$. Using (\ref{conf-g-L-ha}) and an extremal distance argument, we conclude that, as $t\to\infty$, $\ha\xi_L(t)-\Ree\ha g_{L,W}(t,b_1+\pi i)$, $\ha q_{L,m}(t)-\Ree\ha g_{L,W}(t,b_1+\pi i)$, $\Ree\ha g_{L,W}(t,b_2+\pi i)-\ha\xi_L(t)$, and $\Ree\ha g_{L,W}(t,b_2+\pi i)-\ha q_{L,m}(t)$ all tend to $+\infty$. Since $\til g_{L,W}$ has progressive period $(2\pi;2\pi)$, from (\ref{ha-g-L}), $\ha g_{L,W}(t,\cdot)$ has progressive period $(2\ha p; 2\ha p_L(t))$. So $\Ree\ha g_{L,W}(t,b_2+\pi i)-\Ree\ha g_{L,W}(t,b_1+\pi i)=2\ha p_L(t)$. Since $\ha X_{L,m}=\ha\xi_L-\ha q_{L,m}$, we conclude that $2\ha p_L(t)-|\ha X_{{L},m}(t)|\to\infty$ as $t\to\infty$. From Proposition \ref{est-Gamma} (i) and (\ref{pL>p}) we see that $\ln(\ha\Gamma_q(\ha p+t,\ha X_m(t)))$ and $\ln(\ha\Gamma_q(\ha p_{L}(t),\ha X_{{L},m}(t)))$ tend to $0$ as $t\to\infty$. $\Box$

\subsection{Uniformly boundedness} \label{compact-section}
Now we introduce the notation of convergence of domains in \cite{LERW}. We have the following definition and proposition.

\begin{Definition} Suppose $D_n$ is a sequence of plane domains and $D$ is
a plane domain. We say that $(D_n)$ converges to $D$, denoted by $D_n\dto
D$, if for every $z\in D$, $\dist(z,\pa D_n)\to \dist
(z,\pa D)$. This is equivalent to the following: 
\begin{enumerate}
  \item[(i)]  every compact subset of $D$ is contained in all but finitely
many $D_n$'s;
\item [(ii)] for every point $z_0\in\pa D$,
$\dist(z_0,\pa D_n)\to 0$ as $n\to\infty$. 
\end{enumerate}
\label{def-lim}
\end{Definition}

Suppose $D_n\dto D$, and for each $n$, $f_n$ is a complex valued
function on $D_n$, and $f$ is a complex valued function on $D$. We
say that $f_n$ converges to $f$ locally uniformly in $D$, or
$f_n\luto f$ in $D$, if for each compact subset $F$ of $D$, $f_n$
converges to $f$ uniformly on $F$. 

\begin{Proposition} Suppose $D_n\dto D$, $f_n:D_n\conf E_n$ for each $n$, and $f_n\luto f$ in $D$. Then either
$f$ is constant on $D$, or $f$ maps $D$ conformally onto some domain
$E$. And in the latter case, $E_n\dto E$ and $f_n^{-1}\luto f^{-1}$
in $E$. \label{domain convergence}
\end{Proposition}

Fix $m\in\Z$ and $(\rho_1,\rho_2)\in{\cal P}_m$.
Choose $(\rho_1^*,\rho_2^*)\in{\cal P}_m$ such that $\rho_1^*\cup\rho_2^*$ is disjoint from $(\rho_1\cup\rho_2)+2\ha p\Z$, and $\rho_1^*\cup\rho_2^*$ disconnects $\rho_1\cup\rho_2$ from $\ha L$ in $\St_\pi$. Then $\rho_1^*$, $\rho_1$, $\rho_2$,  $\rho_2^*$, and $\rho_1^*+2\ha p$ lie in the order from left to right. Suppose $\rho_j\cap \R=\{a_j\}$, $\rho_j^*\cap\R=\{a_j^*\}$, $\rho_j\cap\R_\pi=\{b_j+\pi i\}$, and $\rho_j^*\cap\R_\pi=\{b_j^*+\pi i\}$, $j=1,2$.  Then we have $a_1^*<a_1<\ha x_0<a_2<a_2^*<a_1^*+2\ha p$ and $b_1^*<b_1<\ha y_m<b_2<b_2^*<b_1^*+2\ha p$.

Let $I_\pi(z)=2\pi-\lin z$ denote the reflection about $\R_\pi$. Let $\Sigma_{\rho_1,\rho_2}$ denote the region in $\St_{2\pi}$ bounded by $\rho_2\cup I_\pi(\rho_2)$ and $(\rho_1\cup I_\pi(\rho_1))+2\ha p$. Fix $\r r\in (b_2^*,b_1^*+2\ha p)$. Then $\r r+\pi i\in\Sigma_H$. Let ${\cal D}_{\rho_1,\rho_2}$ denote the family of simply connected subdomains of $\St_{2\pi}$ which contain $\Sigma_{\rho_1,\rho_2}$, and are symmetric about $I_\pi$. For each $D\in {\cal D}_{\rho_1,\rho_2}$, there is a unique $\r f_D :(\St_{2\pi}; \r r+\pi i)\conf (D;  \r r+\pi i) $ such that $\r f_D'(\r r+\pi i)>0$. Such $\r f_D$ commutes with $I_\pi$. Define a topology on ${\cal D}_{\rho_1,\rho_2}$ such that $D_n\to D_0$ iff $\r f_{D_n}\luto \r f_{D_0}$ in $\St_{2\pi}$.


\begin{Lemma}
 Every sequence   in ${\cal D}_{\rho_1,\rho_2}$ contains a convergent subsequence. \label{comp}
\end{Lemma}
{\bf Proof.} Choose $V$ such that $V:(\St_{2\pi};\r r+\pi i)\conf (\D;0)$. Let $(D_n)$ be a sequence in ${\cal D}_{\rho_1,\rho_2}$. Then $V\circ \r f_{D_n}$, $n\in\N$, is a family of conformal maps from $\St_{2\pi}$ into $\D$. Since this family is uniformly bounded, it contains a subsequence $(V\circ\r f_{D_{n_k}})$ which converges locally uniformly in $\St_{2\pi}$. From Lemma \ref{domain convergence}, this subsequence converges to either a constant function or a conformal map defined on $\St_{2\pi}$. Suppose that the first case happens. Since $V\circ\r f_{D_n}(\r r+\pi i)=0$, the constant is $0$. Then we conclude that $\r f_{D_{n_k}}\luto \r r+\pi i$ in $\St_{2\pi}$, which implies that $\r f_{D_{n_k}}'(\r r+\pi i)\to 0$. Since $\dist (\r r+\pi i,\pa \St_{2\pi})=\pi$, from Koebe's $1/4$ theorem (c.f.\ \cite{Ahl}), we should have $\dist(\r r+\pi i, \pa  {D_{n_k}})\to 0$, which contradicts that
$\dist(\r r+\pi i, \pa  {D_{n_k}})\ge \dist (\r r+\pi i, \R\cup\R_{2\pi}\cup \rho_2\cup(\rho_1+2\ha p))>0$.
 Thus,  $(V\circ \r f_{D_{n_k}})$ converges locally uniformly to a conformal map, which implies that $\r f_{H_{n_k}}$ converges locally uniformly to a conformal map defined on $\St_{2\pi}$, say $\r f$. Since $\r f_{D_{n_k}}$ all map into $\St_{2\pi}$, fix $\r r+\pi i$, have positive derivative at $\r r+\pi i$, and commute with $I_\pi$, $\r f$ should also satisfy these properties. Let $D_0=\r f(\St_{2\pi})$. Then $D_0$ is a simply connected subdomain of $\St_{2\pi}$, contains $\r r+\pi i$, and is symmetric about $\R_\pi$. We suffice to show that $D_0\supset \Sigma_{\rho_1,\rho_2}$ because if this is true, then $D_0\in {\cal D}_{\rho_1,\rho_2}$ and $\r f=\r f_{D_0}$, which implies that $D_{n_k}\to D_0$. Suppose $D_0\not\supset \Sigma_{\rho_1,\rho_2}$. Since $\r r+\pi i\in  D_0$, and $\Sigma_{\rho_1,\rho_2}$ is connected, there exists $z_0\in \Sigma_{\rho_1,\rho_2}\cap D_0$. From Lemma \ref{domain convergence} we have $ {D_{n_k}}\dto D_0$. From Definition \ref{def-lim} (ii), we see that $\dist(z_0,\pa D_{n_k})\to 0$, which contradicts that $z_0\in\Sigma_{\rho_1,\rho_2}\subset D_{n_k}$ for each $k$. This finishes the proof. $\Box$

\vskip 4mm

Let $I_0(z)=\lin z$ denote the reflection about $\R$. For $D\in {\cal D}_{\rho_1,\rho_2}$, let
$$D^\pm=D\cup I_0(D)\cup (a_2,a_1+2\ha p).$$
Then $D^\pm$ is a simply connected subdomain of $\St_{2\pi}^\pm:=\{-2\pi<\Imm z<2\pi\}$, and is symmetric about $\R$. Let $\r g_D=\r f_D^{-1}:D\conf \St_{2\pi}$.  From Schwarz reflection principle, $\r g_D$ extends to a conformal map $\r g_D^\pm$ from $D^\pm$ into $\St_{2\pi}^{\pm}$, which commutes with $I_0$.

\begin{Lemma}
  If $D_n\to D_0$, then $D_n^\pm\dto D_0^\pm$ and $\r g^\pm_{D_n}\luto \r g_{D_0}^\pm$ in $D_0^\pm$. \label{comp2}
\end{Lemma}
{\bf Proof.} From Lemma \ref{domain convergence} we have $D_n\dto D_0$ and $\r g_{D_n}\luto \r g_{D_0}$ in $D_0$. Then we easily see that $D_n^\pm\dto D_0^\pm$. Let $(D_{n_k})$ be a subsequence of $(D_n)$. Choose $V:\St_{2\pi}^\pm\conf\D$. Then $(V\circ \r g_{D_{n_k}}^\pm)$ is uniformly bounded family, which contains a subsequence $(V\circ \r g_{D_{n_{k_l}}}^\pm)$ that converges locally uniformly to some function $G$ in $D_0^\pm$. Since $\r g_{D_{n_{k_l}}}^\pm\luto g_{D_0}^\pm$ in $D_0$, we see that $G$ is the analytic extension of $V\circ g_{D_0}$. Thus, $G=V\circ g_{D_0}^\pm$. So we conclude that $\r g_{D_{n_{k_l}}}^\pm\luto g_{D_0}^\pm$ in $D_0^\pm$. The proof is now finished because every subsequence of $(\r g^\pm_{D_n})$ contains a subsequence which converges to $\r g_{D_0}^\pm$ locally uniformly in $D_0^\pm$. $\Box$

\vskip 4mm

For each $D\in {\cal D}_{\rho_1,\rho_2}$, let $D({\ha L})$ be the connected component of $D\sem(\ha L\cup I_\pi(\ha L))$ that contains $\r r+\pi i$. Then $D({\ha L})$ is a simply connected subdomain of $\St_{2\pi}$, and is symmetric about $\R_\pi$. There is a unique $\r g_{D,\ha L}$ such that $\r g_{D,\ha L}: (D({ \ha L}); \r r+\pi i)\conf (\St_{2\pi}; \r r+\pi i)$ and $\r g_{D,\ha L}'(\r +\pi i)>0$.
Let  $$D^\pm({\ha L})=D({ \ha L})\cup I_0(D({ \ha L}))\cup ( (a_2,a_1+2\ha p)\sem \lin{\ha L}).$$
 Then $\r g_{D,\ha L}$ extends to a conformal map $\r g^\pm_{D,\ha L}$ from $D^\pm({\ha L})$ into $\St_{2\pi}^\pm$, which commutes with $I_0$. We easily see that $D_n\dto D_0$ iff $D_n^\pm(\ha L)\dto D_0^\pm (\ha L)$. Using some subsequence argument we can derive the following lemma.

\begin{Lemma}
  If $D_n\to D_0$, then $D_n^\pm(\ha L)\dto D_0^\pm(\ha L)$ and $\r g^\pm_{D_n,\ha L}\luto \r g_{D_0,\ha L}^\pm$ in $D_0^\pm(\ha L)$. \label{comp3}
\end{Lemma}

For each $D\in {\cal D}_{\rho_1,\rho_2}$, $\lin{\rho_1^*}$ and $\lin{\rho_2^*}$ are compact subsets of $D^\pm$ and $D^\pm({\ha L})$. From Lemma \ref{comp}, Lemma \ref{comp2}, and Lemma \ref{comp3} we conclude that, there is a constant $C>0$ which depends only on $\rho_1,\rho_2,\rho_1^*,\rho_2^*,\ha L,\r r$ such that, for any $D\in {\cal D}_{\rho_1,\rho_2}$ and $z\in\rho_1^*\cup\rho_2^*$, the following quantities:
$$|\r g_D(z)-z|, \quad |\r g_D'(z)|, \quad |1/\r g_D'(z)|, \quad |\r g_{D,\ha L}(z)-z|,\quad |\r g_{D,\ha L}'(z)|,\quad |1/\r g_{D,\ha L}'(z)|,$$
are all bounded above by $C$.

Fix $t\in[0,\ha T_{\rho_1,\rho_2})$. Let $D=\St_{2\pi}\sem(\ha\beta((0,t])+2\ha p\Z)\sem I_\pi(\ha\beta((0,t])+2\ha p\Z)$. Then $D\in {\cal D}_{\rho_1,\rho_2}$. We have $\r g_D: (\St_\pi\sem (\ha\beta((0,t])+2\ha p\Z);\R_\pi)\conf (\St_\pi;\R_\pi)$. Let $h_1=\ha g(t,\r r+ \pi i)-(\r r+\pi i)\in\R$. Since $\r g_D$ fixes $\r r+\pi i$, from  (\ref{conf-g-ha}) we have $\r g_D=\ha g(t,\cdot)-h_1$. Similarly, using (\ref{conf-g-L-ha}) we conclude that $\r g_{D,\ha L}=\ha g_{L,W}(t,\cdot)-h_2$ for some $h_2\in\R$. Thus, for any $z\in\rho_1^*\cup\rho_2^*$, the following quantities:
$$|\ha g(t,z)-z-h_1|, \quad |\ha g'(t,z)|, \quad |1/\ha g'(t,z)|, \quad |\ha g_{L,W}(t,z)-z-h_2|,\quad |\ha g_{L,W}'(t,z)|,\quad |1/\ha g_{L,W}'(t,z)|,$$
are all bounded above by the $C$ in the last paragraph. Let $h=h_2-h_1$ and $C'=\max\{2C,C^2\}$. From (\ref{W-ha}), we find that,
\BGE |\ha W(t,z)-z-h|\le C',\quad 1/C'\le |\ha W'(t,z)|\le C',\quad z\in\ha g(t,\rho_1^*)\cup\ha g(t,\rho_2^*).\label{haW-est}\EDE

\no{\bf Proof of Proposition \ref{N-to-1-ha} (ii).} From (\ref{N-bd}), we suffice to show that (ii) holds if $\ln(\ha N_m(t))$ is replaced by $\ln(\ha A_1(t))$, $\ln(\ha A_{I,m}(t))$, $\ln(\ha\Gamma_q(\ha p+t,\ha X_m(t)))$, $\ln(\ha\Gamma_q(\ha p_{L}(t),\ha X_{{L},m}(t)))$, $\ha X_{{L},m}(t)-\ha X_m(t)$, or $\ha \RA(\ha p_{L}(t))-\ha \RA(\ha p+t)$, respectively. From (\ref{RA-rA-ha}) and (\ref{pL>p}) we see that  $\ha \RA(\ha p_{L}(t))$ and $\ha \RA(\ha p+t)$ are both positive and bounded above by $\ha\RA(\ha p)$, which is a uniform constant. So $\ha \RA(\ha p_{L}(t))-\ha \RA(\ha p+t)$ is uniformly bounded.

Fix $(\rho_1,\rho_2)\in {\cal P}_m$ and $t\in [0,\ha T_{\rho_1,\rho_2})$. Then (\ref{haW-est}) holds.
From Schwarz reflection principle, $\ha W(t,\cdot)$ extends conformally to a conformal map on
$\Sigma:=\C\sem (\lin{\ha L_t}\cup I_0(\ha L_t)+2\pi i \Z)$,
and the extended map commutes with both $I_0$ and $I_\pi$. Thus, $\ha W(t,\cdot)$ has progressive period $(2\pi i;2\pi i)$. So $\ha W'(t,\cdot)$, $1/\ha W'(t,\cdot)$, and $\ha W(t,\cdot)-\cdot$ are all analytic functions with period $2\pi i$. Let $$\rho_{j,t}^*=(\ha g(t, \rho_1^*)\cup I_0(\ha g(t, \rho_1^*)))+2\pi i\Z,\quad j=1,2.$$ Then $\rho_{1,t}^*$ and $\rho_{2,t}^*$ are two disjoint simple curves with period $2\pi i$, which lie inside $\Sigma$, and (\ref{haW-est}) holds for any $z\in\rho_{1,t}^*\cup\rho_{2,t}^*$.
Since $\ha\xi(t)$ and $\ha q_m(t)+\pi i$ lie inside the region bounded by $\rho_{1,t}^*$ and $\rho_{2,t}^*$,  from Maximum Principle, (\ref{ha-A-def}), and (\ref{W(q)-ha-m})  we have $$|\ha\xi_L(t)-\ha\xi(t)-h|,|\ha q_{L,m}(t)-\ha q_m(t)-h|\le C',\quad
1/C'\le \ha A_1(t), \ha A_{I,m}(t)\le C'.$$
Since $\ha X_m=\ha\xi-\ha q_m$ and $\ha X_{L,m}=\ha\xi_L-\ha q_{L,m}$, we have $|\ha X_{L,m}(t)-\ha X_m(t)|\le 2 C'$.
Thus, the lemma holds if $\ln(\ha N_m(t))$ is replaced by $\ln(\ha A_1(t))$, $\ln(\ha A_{I,m}(t))$, or $\ha X_{{L},m}(t)-\ha X_m(t)$.

We know that $\rho_1^*$, $\rho_2^*$, and $\rho_1^*+2\ha p$ are pairwise disjoint, and lie in the order from left to right.
Since $\til g(\check t,\cdot)$ has progressive period $(2\pi;2\pi)$, from (\ref{ha-g-L}), $\ha g(t,\cdot)$ has progressive period $(2\ha p;2(\ha p+t))$. Thus, $\ha g(t,\rho_1^*)$, $\ha g(t,\rho_2^*)$, and $\ha g(t,\rho_1^*)+2(\ha p+t)$ are pairwise disjoint, and lie in the order from left to right. Since $\ha \xi(t)$ and $\ha q_m(t)+\pi i$ are bounded by $\ha g(t,\rho_1^*)$ and $\ha g(t,\rho_2^*)$ in $\St_\pi$, they are also bounded by $\ha g(t,\rho_1^*)$ and $\ha g(t,\rho_1^*)+2(\ha p+t)$ in $\St_\pi$. Thus,   $|\ha X_m(t)|=|\ha\xi(t)-\ha q_m(t)|$ is bounded above by $2(\ha p+t)+\diam(\ha g(t,\rho_1^*))$. Since $|\ha g'(t,z)|\le C$ on $\rho_1^*$,   $\diam(\ha g(t,\rho_1^*)\le C\diam(\rho_1^*)$. Thus, $|\ha X_m(t)|-2(\ha p+t)$ is bounded above by a uniform constant. From Proposition \ref{est-Gamma} (ii) we see that the lemma holds if $\ln(\ha N_m(t))$ is replaced by $\ln(\ha\Gamma_q(\ha p+t,\ha X_m(t)))$. Similarly, $|\ha X_{L,m}(t)|-2\ha p_L(t)$ is bounded above by a uniform constant, which implies that the lemma holds if $\ln(\ha N_m(t))$ is replaced by $\ln(\ha\Gamma_q(\ha p_{L}(t),\ha X_{{L},m}(t)))$. $\Box$

\section{Restriction}\label{Sec7}
\subsection{Brownian loop measure} \label{Section-B-L}

\begin{Lemma}
  Let $p_0>0$ and $L_0$ be a  hull in $\A_{p_0}$ w.r.t.\ $\TT_{p_0}$. Let $\til L_0=(e^i)^{-1}(L_0)$. Suppose that $p_1=\modd(\A_{p_0}\sem L)\in(0,p_0)$,  $W_0:(\A_{p_0}\sem L_0;\TT_{p_0})\conf(\A_{p_1};\R_{p_1})$,
and  $\til W_0:(\St_{p_0}\sem \til L_0;\R_{p_0})\conf(\St_{p_1};\R_{p_1})$,
 and $e^i\circ \til W_0=W_0\circ e^i$.
Let $x\in\R$ be such that $\dist(e^{ix},L_0)>0$. 
Let $S\til W_0(x_0)$ denote the Schwarz derivative of $\til W_0$ at $x_0$. Let $\mu_{e^{ix_0}}$ denote the Brownian bubble measure in $\A_{p-t}$ rooted at $e^{ix_0}$. Let ${\cal E}_{L_0}$ denote the set of curves that intersect $L_0$. Then
  $$\mu_{e^{ix_0}}[{\cal E}_{L_0}]=-\frac 16 S\til W_0(x_0)+\frac 12\til W_0'(x_0)^2(\rA(p_1)+\frac 1{p_1})-\frac 12(\rA(p_0)+\frac 1{p_0}).$$
  \label{BL0}
\end{Lemma}
{\bf Proof.} Let $z_0\in\St_{p_0}$. The bubble measure $\mu_{e^{ix_0}}$ equals $\lim_{z_0\to x_0}\frac{\PP_{{z_0};{x_0}}} {|z_0-x_0|^2}$, where $\PP_{{z_0};{x_0}}$ is the distribution of a planar Brownian motion started from $e^{iz_0}$ conditioned to exit $\A_{p_0}$ from $e^{ix_0}$. Choose $x_1,x_2\in\R$ such that $x_1<x_0<x_2<x_1+2\pi$. Then $\PP_{{z_0};{x_0}}$ equals the limit of $\PP_{{z_0};(x_1,x_2)}$ as $x_1,x_2\to x_0$, where $\PP_{{z_0};(x_1,x_2)}:=\PP_{{z_0}}[\cdot|{\cal E}_{x_1,x_2}]$, $\PP_{{z_0}}$ is the distribution of a planar Brownian motion started from $e^{iz_0}$, and ${\cal E}_{x_1,x_2}$ denotes the event that the curve ends at the arc $e^i((x_1,x_2))$.

Since the Poisson kernel function in $\A_{p_0}$ with the pole at $e^{ix}\in\TT$ is
$z\mapsto\frac 1{2\pi}(\Ree\SA(p_0,z/e^{ix})+\frac{\ln|z|}{p_0})$, we get
\BGE \PP_{{z_0}}[{\cal E}_{x_1,x_2}]=-\frac 1{2\pi} \int_{x_1}^{x_2}\Imm (\HA(p_0,z_0-x)+\frac{z_0}{p_0})dx.\label{PEx1x2}\EDE
From conformal invariance of planar Brownian motions, $\PP_{{z_0}}[{\cal E}_{x_1,x_2}\sem {\cal E}_{L_0}]$ is equal to the probability of a planar Brownian motion started from $W_0(e^{iz_0})=e^i(\til W_0(z_0))$ hits $\pa \A_{p_1}$ at the arc $W_0(e^i((x_1,x_2)))=e^i((\til W_0(x_1),\til W_0(x_2)))$. From (\ref{PEx1x2}) and change of variables, we get
$$\PP_{{z_0}}[{\cal E}_{x_1,x_2}\sem{\cal E}_{L_0}]=-\frac 1{2\pi} \int_{x_1}^{x_2}\Imm (\HA(p_1,\til W_0(z_0)-\til W_0(x))+\frac{\til W_0(z_0)}{p_1}) \til W_0'(x)dx.$$ 

Then we get an expression for $\PP_{z_0;x_1,x_2}[{\cal E}_{L_0}]=\PP_{z_0}[{\cal E}_{L_0}|{\cal E}_{x_1,x_2}$.
Letting $x_1,x_2\to x_0$, we get
$$ \PP_{z_0;x_0}[{\cal E}_{L_0}]=1-\frac{\til W_0'(x_0) \Imm (\HA(p_1,\til W_0(z_0)-\til W_0(x_0))+\frac{\til W_0(z_0)-\til W_0(x_0)}{p_1}) }{\Imm (\HA(p_0,z_0-x_0)+\frac{z_0-x_0}{p_0})}.$$ 
Finally we compute $\lim_{z_0\to x_0}\frac{\PP_{{z_0};{x_0}}[{\cal E}_{L_0}]} {|z_0-x_0|^2}$.
The proof is completed by some tedious but straightforward computation involving power series expansions. $\Box$


\begin{Lemma}
For the $U(t)$ defined in (\ref{U(t)}), we have  $\mu_{\mbox{loop}}[{\cal L}_{L,t}]=U(t)$, $0\le t\le T$, where ${\cal L}_{L,t}$ denotes the set of loops in $\A_p$ that intersect both $L$ and  $\beta((0,t))$. \label{BL}
\end{Lemma}
{\bf Proof.} For $0\le t<T$, let $\mu_t$ denote the Brownian bubble measure in $\A_{p-t}$ rooted at $e^{i\xi(t)}$. The argument in \cite{loop} shows that $\mu_{\mbox{loop}}[{\cal L}_{L,t}]=\int_0^t \mu_s[\{\cdot\cap L_s\ne\emptyset\}]ds$ for $0\le t\le T$.
From (\ref{conf-W-til}), (\ref{A1I}) and the previous lemma, we have
$$\mu_s[\{\cdot\cap L_s\ne\emptyset\}]=-\frac 16 A_S(s)+\frac 12A_1(s)^2(\rA(p_L-v(s))+\frac 1{p_L-v(s)})-\frac 12(\rA(p-s)+\frac 1{p-s}).$$
The proof can now be completed by integrating the right-hand side of this formula from $0$ to $t$ and using (\ref{v'*}) and that $v(0)=0$. $\Box$

\begin{Lemma} Let $m\in\Z$.
\begin{enumerate}
  \item [(i)]  On the event ${\cal E}_m\cap\{\beta\cap L=\emptyset\}$, $U(p)$ is finite.
  \item [(ii)] For any $(\rho_1,\rho_2)\in \til{\cal P}_m$, $U(t)$ is uniformly bounded on $[0,\til T_{\rho_1,\rho_2})$.
\end{enumerate}
  \label{R-to-end}
\end{Lemma}
{\bf Proof.} From \cite{loop}, if two sets in $\C$ have positive distance from each other, then the Brownian loop measure of the loops that intersect both of them is finite.  (i) If ${\cal E}_m$ occurs and $\beta\cap L=\emptyset$, then  $\dist(L,\beta((0,p)))>0$. From Lemma \ref{BL} and the above observation, $U(p)=\mu_{\mbox{loop}}[{\cal L}_{L,p}]$ is finite. (ii)  Let ${\cal L}_{L,\rho_1,\rho_2}$ denote the set of loops in $\A_p$ that intersect both $L$ and ${\rho_1\cup\rho_2}$. Since $\dist(L,{\rho_1\cup\rho_2})>0$, we have $\mu_{\mbox{loop}}[{\cal L}_{L,\rho_1,\rho_2}]<\infty$. If $t< T_{\rho_1,\rho_2}$, then $\rho_1\cup\rho_2$ disconnects $\beta((0,t])$ from $L$, which means that a loop in $\A_p$ that intersects both $L$ and $\beta((0,t])$ must also intersect $\rho_1\cup\rho_2$. So ${\cal L}_{L,t}\subset {\cal L}_{L,\rho_1,\rho_2}$. Thus, $U(t)$, $0\le t<T$, is bounded above by $\mu_{\mbox{loop}}[{\cal L}_{L,\rho_1,\rho_2}]$. $\Box$

\subsection{Radon-Nikodym derivatives}  \label{Section-RN}
Let $s\in\R$ and $m\in\Z$. Consider the following two SDEs:
\BGE d\xi(t)=\sqrt\kappa dB(t)+\Big(3-\frac\kappa 2\Big)\frac{A_2(t)}{A_1(t)}dt+{A_1(t)} {\Lambda_{\langle s \rangle}(p_L-v(t),X_{L,0}(t))}dt,\quad 0\le t<T;\label{xi-K}\EDE
\BGE d\xi(t)=\sqrt\kappa dB(t)+\Big(3-\frac\kappa 2\Big)\frac{A_2(t)}{A_1(t)}dt+{A_1(t)} {\Lambda_0(p_L-v(t),X_{L,m}(t))}dt,\quad 0\le t<T.\label{xi-K-m}\EDE
Let the distribution of $(\xi(t), 0\le t<T)$ be denoted by  $\mu_{L,\langle s\rangle}$ or $\mu_{L,m}$, respectively, if $(\xi(t))$, $0\le t<T$, is the maximal solution of (\ref{xi-K}) or (\ref{xi-K-m}), respectively, and   $\xi(0)=x_0$.

Suppose  that $(\xi)$ has distribution $\mu_{L,m}$. From (\ref{W(t)}),  (\ref{qLm}),  (\ref{X}) and  (\ref{dxiK}),  we get
$$d\xi_L(t) =A_1(t)\sqrt\kappa dB(t)+A_1(t)^2  {\Lambda_0(p_L-v(t),\xi_L(t)-\Ree \til g_L(t,\til W_L(y_m+ p i)))}dt,\quad 0\le t<T.$$
Since  $\xi_L(t)=\eta_L(v(t))$ and $\til g_L(t,\cdot)=\til h_L(v(t),\cdot)$, from (\ref{v'*}) and (\ref{A1I})
we conclude that there is another Brownian motion $B_{v}(t)$ such that
$$d\eta_L(t)=\sqrt\kappa dB_{v}(t)+{\Lambda_0(p_L-t,\eta_L(t)-\Ree\til h_L(t,\til W_L(y_m+p i)))}dt,\quad 0\le t<v(T).$$
Recall that $\til h_L$ and $\til\gamma_L$ are the covering annulus Loewner maps and trace of modulus $p_L$ driven by $\eta_L$. Thus, $\til\gamma_L(t)$, $0\le t<v(T)$, is a covering annulus SLE$(\kappa;\Lambda_0)$ trace in $\St_{p_L}$ started from
$\til W_L(\xi(0))$ with marked point $\til W_L(y_m+p i)$, stopped at $v(T)$.

There are two possibilities. Case 1: $v(T)=p_L$. Then $\til\gamma_L(t)$, $0\le t<v(T)$ is a complete covering annulus SLE$(\kappa;\Lambda_0)$ trace. From the last paragraph of Section \ref{section-limit} we know that a.s.\ $\lim_{t\to v(T)^-} \til\gamma_L(t)= \til W_L(y_m+p i)$. Since $\til\gamma_L(t)=\til\beta_L(v(t))=\til W_L(\beta(v(t)))$, we have $T=p$ and $\lim_{t\to T^-}\til\beta(t)=y_m+ p i$, which means that the event ${\cal E}_m$ occurs. Case 2: $v(T)<p_L$. Then $\lim_{t\to v(T)^-} \til\gamma_L(t)$ exists and lie in $\St_{p_L}$, which implies that $\lim_{t\to T^-} \til\beta(t)$ exists and lie in $\St_p\sem \til L$. This means that the solution $\xi(t)$, $0\le t<T$, can be further extended, which is a contradiction. So only Case 1 can happen, which implies that $\mu_{L,m}(\{T=p\}\cap {\cal E}_m)=1$. 

Similarly, if $(\xi(t))$ has the distribution $\mu_{L,\langle s\rangle}$, then a.s.\ $v(T)=p_L$, $\til\gamma_L(t)$, $0\le t<v(T)$, is a complete covering annulus SLE$(\kappa;\Lambda_{\langle s\rangle})$ trace in $\St_{p_L}$ started from $\til W_L(\xi(0))$ with marked point $\til W_L(y_0+p i)$, and $\lim_{t\to v(T)^-} \til\gamma_L(t)$ exists and belongs to $y_0+ pi+ 2\pi\Z$. Thus, $\mu_{L,\langle s\rangle}(\{T=p\})=1$ and $\mu_{L,\langle s\rangle}(\bigcup_{m\in\Z} {\cal E}_m)=1$. Since $X_{L,m}(0)=\til W_L(\xi(0))-\Ree \til W_L(y_m+p i)$,
from (\ref{dmu-m/dmu}) we have \BGE \frac{d\mu_{L,m}}{d\mu_{L,\langle s\rangle}}=e^{\frac{2\pi}\kappa ms}\, \frac{ \Gamma_0(p_L,X_{L,m}(0))}{\Gamma_{\langle s\rangle}(p_L,X_{L,0}(0))}\,{\bf 1}_{{\cal E}_m}.\label{dmu-m/dmu-L}\EDE

Suppose $(\xi(t))$ has the distribution $\mu_{L,\langle s\rangle}$. Since $\gamma_L$ is the trace driven by $\eta_L$, the above argument shows that, $\gamma_L(t)$, $0\le t<v(T)$, is a complete annulus SLE$(\kappa;\Lambda_{\langle s\rangle})$ trace in $\A_{p_L}$ started from $e^i\circ \til W_L(\xi(0))=W_L(e^{ix_0})$ with marked point $e^i\circ \til W_L(y_0+p i)=W_L(e^{iy_0-p})$. Since $W_L:(\A_p\sem L;\beta(v^{-1}(t)))\conf (\A_{p_L};\gamma_L(t))$,   we see that $\beta(v^{-1}(t))$, $0\le t<v(T)$, is an annulus SLE$(\kappa;\Lambda_{\langle s\rangle})$ trace in $\A_p\sem L$ started from $e^{ix_0}$ with marked point $e^{iy_0-p}$.

The process $(M_m(t))$ defined earlier will be used to derive the Radon-Nikodym derivative between the $\mu_{L,m}$ defined here and the $\mu_m$ defined as the distribution of the solution of (\ref{xi-p-m}). Suppose that $(\xi(t))$ has distribution $\mu_m$. Then $\xi(t)$, $0\le t<p$, solves the SDE: \BGE d\xi(t)=\sqrt\kappa dB(t)+\Lambda_0(p-t,X_m(t))dt,\quad 0\le t<p,\quad \xi(0)=x_0,\label{xi-p-m-X}\EDE
From (\ref{dM/M}) we see that $M_m(t)$, $0\le t<T$, is a local martingale under $\mu_m$.

Let $(\rho_1,\rho_2)\in\til {\cal P}_m$. From Propostion \ref{N-to-1} (ii),  Lemma \ref{R-to-end} (ii), and $M_m=N_m\exp(\cc U)$, we see that $M_m(t)$ is uniformly bounded on $[0,T_{\rho_1,\rho_2})$. Thus, $M_m(t\wedge T_{\rho_1,\rho_2})$ is a bounded martingale, and we have
$\EE_{\mu_m}[M_m(T_{\rho_1,\rho_2})]=M_m(0)$. If we now change the distribution of $(\xi(t))$ from $\mu_m$ to a new probability measure $\nu$ defined by $d\nu/d\mu_m=M_m(T_{ \rho_1, \rho_2})/M_m(0)$, then from Girsanov's Theorem and (\ref{dM/M}) we see that the current $\xi(t)$ satisfies SDE (\ref{xi-K-m}) for $0\le t<T_{ \rho_1, \rho_2}$. Thus, on the event $\{T_{ \rho_1, \rho_2}=p\}$, $\mu_{L,m}\ll\mu_m$, and the Radon-Nikodym derivative between the two measures restricted to the event $\{T_{ \rho_1, \rho_2}=p\}$ is $M_m(p)/M_m(0)$. From Proposition \ref{N-to-1}(i), Lemma \ref{R-to-end} (i), (\ref{T=p}) and $M_m= N_m\exp(\cc U)$, we see that $M_m(p)=C_{p,L}\exp(\cc U(p))$. So
\BGE d\mu_{L,m}/d\mu_m=C_{p,L}\exp(\cc U(p))/M_m(0)\quad\mbox{on}\quad\{T_{ \rho_1, \rho_2}=p\}.\label{dmuL/dmum-rho}\EDE

Suppose ${\cal E}_m$ occurs and $T=p$. Then $\til \beta \cap \til L=\emptyset$. Since $\til \beta$ starts from $x_0$, we can find $( \rho_1, \rho_2)\in\til{\cal P}_m$ such that $\til\beta\cap(\rho_1\cup\rho_2)=\emptyset$, which implies that $ T_{\rho_1,\rho_2}=p$. Thus,
\BGE {\cal E}_m\cap\{T=p\}\subset \bigcup_{( \rho_1, \rho_2)\in\til{\cal P}_m} \{ T_{ \rho_1, \rho_2}=p\}.\label{union}\EDE
Since $\mu_m({\cal E}_m)=1$ and ${\cal P}_m$ is countable, from (\ref{dmuL/dmum-rho}) and (\ref{union}) we see that $d\mu_{L,m}/d\mu_m=C_{p,L}\exp(\cc U(p))/M_m(0)$ on $\{T=p\}$. Since $\mu_{L,m}(\{T=p\})=1$, from (\ref{T=p}) and Lemma \ref{BL} we know that
\BGE \frac{d\mu_{L,m}}{d\mu_m}=\frac{C_{p,L}}{M_m(0)}{\bf 1}_{\{\til\beta\cap  \til L=\emptyset\}}  \exp(\cc \mu_{\mbox{loop}}[{\cal L}_{L,p}]),\label{dmuLm/dmum}\EDE
where  ${\cal L}_{L,p}$ is the set of loops in $\A_p$ that intersect both $L$ and  $\beta((0,p))$.

Let $s\in\R$. Now we compare $\mu_{\langle s\rangle}$ with $\mu_{L,\langle s\rangle}$. Define
$$ Y_{\langle s\rangle}(t)=\Gamma_{\langle s\rangle}(p-t,X_0(t)),\quad Y_{L,{\langle s\rangle}}(t)=\Gamma_{\langle s\rangle}(p_L-v(t),X_{L,0}(t)).$$ 
Define $M_{\langle s\rangle}$ using (\ref{Mm}) with $Y_m$ and $Y_{L,m}$ replaced by $Y_{\langle s\rangle}$ and $Y_{L,\langle s\rangle}$, respectively, and $A_{I,m}$ replaced by $A_{I,0}$.

Since $\til g(t,\cdot)$ has progressive period $(2\pi;2\pi)$, from (\ref{qm}) we have $q_m(t)=q_0(t)+2m\pi$.
Since $\til W(t,\cdot)$ has progressive period $(2\pi;2\pi)$, from (\ref{A1I}) we have $A_{I,m}=A_{I,0}$. Thus,
$$\frac{M_m(0)}{M_{\langle s\rangle}(0)}=\frac{ \Gamma_0(p_L,X_{L,m}(0))/ \Gamma_0(p,X_m(0))}{\Gamma_{\langle s\rangle}(p_L,X_{L,0}(0))/\Gamma_{\langle s\rangle}(p,X_0(0))}.$$
Since $X_m(0)=x_0-y_m$, from (\ref{dmu-m/dmu}), (\ref{dmu-m/dmu-L}),   (\ref{dmuLm/dmum}) and the above formula, we get
 $$\frac{d\mu_{L,\langle s\rangle}}{d\mu_{\langle s\rangle}}=\frac{C_{p,L}}{M_{\langle s\rangle}(0)} {\bf 1}_{\{\beta\cap  L=\emptyset\}} \exp (\cc \mu_{\mbox{loop}}[{\cal L}_{L,p}]).$$ 
Recall that when $(\xi(t))$ has distribution $\mu_{\langle s\rangle}$, $\beta(t)$, $0\le t<p$, is an annulus SLE$(\kappa;\Lambda_{\langle s\rangle})$ trace in $\A_p$ started from $z_0=e^{ix_0}$ with marked point $w_0=e^{iy_0-p}$. When  $(\xi(t))$ has distribution $\mu_{L,\langle s\rangle}$, a time change of $\beta$: $\beta(v^{-1}(t))$, $0\le t<v^{-1}(T)$, is an annulus SLE$(\kappa;\Lambda_{\langle s\rangle})$ trace in $\A_p\sem L$ started from $z_0$ with marked point $w_0$. So we finish the proof of Theorem \ref{main-ann}.

\section{Other Results}\label{Sec8}
\subsection{Restriction in a simply connected subdomain}
We now give a sketch of the proof of Theorem \ref{Main-strip}. Let $p>0$, $\kappa\in(0,4]$, $s\in\R$ $z_0\in\TT$, $w_0\in\TT_p$, and the set $L$ be as in Theorem \ref{main-ann}. Choose $x_0,y_0\in\R$ such that $z_0=e^{ix_0}$ and $w_0=e^{iy_0-p}$. Let $y_m=y_0+2m\pi$, $m\in\Z$. Let $\til L=(e^i)^{-1}(L)$. Then $\St_p\sem \til L$ is a disjoint union of simply connected domains $\til D_m$, $m\in\Z$, such that $\til D_m=\til D_0+2m\pi$ for $m\in\Z$. We label one of the domains $\til D_0$ such that  $x_0\in\pa\til D_0$. There is a unique $m_0\in\Z$ such that $y_{m_0}+ p i\in\pa \til D_0$. We have $e^i:\til D_0\conf \A_p\sem L$.
Let $J_0$ be the component of $\TT_p\sem \lin{L}$ that contains $w_0$. We may find $W_L$ such that $W_L:(\A_p\sem L;J_0)\conf (\St_\pi;\R_\pi)$. Let $\til W_L=W_L\circ e^i$, and $\til J_0$ be a component of $\R_p\sem \lin{\til L}$ that contains $y_{m_0}+ p i$. Then $\til W_L:(\til D_0;\til J_0)\conf (\St_\pi;\R_\pi)$.

Let $\xi(t)$, $g(t,\cdot)$, $\til g(t,\cdot)$, $\beta(t)$, $\til\beta(t)$,  $0\le t<p$, and $T\in (0,p]$ be as in Section \ref{Section-SDE}. Now we define $\til\beta_L(t)=W_L(\beta(t))=\til W_L(\til\beta(t))$, $0\le t<T$. Then $\til\beta_L$ is a simple curve with $\til\beta(0)\in\R$ and $\til\beta((0,p))\subset\St_\pi$. Let $v(t)$ be the capacity of $\til\beta_L((0,t])$ in $\St_\pi$ w.r.t.\ $\R_\pi$ for $0\le t<T$. Let $S=\sup v([0,T))$, and $\til\gamma_L(t)=\til\beta_L(v^{-1}(t))$, $0\le t<S$. Then $\til\gamma_L$ is the strip Loewner trace driven by some $\eta_L\in C([0,S))$.

Let $\til h_L(t,\cdot)$, $0\le t<S$, be the strip Loewner maps  driven by $\eta_L$. Define $\xi_L(t)=\eta_L(v(t))$ and $\til g_L(t,\cdot)=\til h_L(v(t),\cdot)$. Define $\til g_{L,W}(t,\cdot)$ and $\til W(t,\cdot)$ using (\ref{W(t)}). Then (\ref{conf-g-L-til}) and (\ref{conf-W-til}) hold with $p_L-v(t)$ replaced by $\pi$. From (\ref{trace-strip}) we see that (\ref{trace-L}) and (\ref{W(xi)}) hold.

For $m\in\Z$, define $q_m(t)$ and $q_{L,m}(t)$ using (\ref{qm}) and (\ref{qLm}) with $p_L-v(t)$ replaced by $\pi$. Define $A_j(t)$ and $A_{I,m}(t)$ using (\ref{A1I}).  Define $X_m(t)$ and $X_{L,m}(t)$ using (\ref{X}). A standard argument shows that (\ref{v'*}) holds here. So (\ref{patg-L}) holds with $\HA(p_L-v(t),\cdot)$ replaced by $\coth_2$.
 Now (\ref{dq(t)}) and  (\ref{gI'}) still hold here. From (\ref{ODE-HA-I-strip}) and (\ref{deriv2-strip}) we see that (\ref{dqK}), (\ref{gKI'}) and (\ref{dAI}) hold here with $\HA_I(p_L-v(t),\cdot)$ replaced by $\tanh_2$.


By differentiating  $\til W(t,\cdot)\circ \til g(t,z)=\til g_{L,W}(t,z)$ w.r.t.\ $t$ and $z$, and letting $w=\til g(t,z)\to\xi(t)$, we conclude that (\ref{-3*}) holds here, and (\ref{patW'}) holds with $\rA(p_L-v(t))$ replaced by $\frac 16$, which comes from the power series expansion: $\coth_2(z)=\frac 2z+\frac z6 +O({z^2})$ when $z$ is near $0$.  
Then (\ref{dxiK}) and  (\ref{dXt}) still hold here;  (\ref{dXK}) holds with $\HA_I(p_L-v(t),\cdot)$ replaced by $\tanh_2$; and (\ref{dA1}) should be modified with $\frac 16$ in place of $\rA(p_L-v(t))$.

Define $Y_m(t)$ using (\ref{Y=X-m}), but define $Y_{L,m}(t):=\ha\Gamma_\infty(v(t),X_{L,m}(t))$.
Since $\Gamma_0$ solves (\ref{PDE-Gamma}) and $\ha\Gamma_\infty$ solves (\ref{PDE-Gamma-infty}), using (\ref{v'*}), (\ref{dXt}), and the modified (\ref{dAI}) and (\ref{dXK}) we find that (\ref{dY-m}) still holds, and (\ref{dYK-m}) holds with  $\HA_I(p_L-v(t),\cdot)$  and $\Lambda_0(p_L-v(t),\cdot)$ replaced by replaced by $\tanh_2$ and $\kappa {\ha\Gamma_\infty'}/{\ha\Gamma_\infty}=(\frac\kappa 2-3)\tanh_2$, respectively.

Define $M_m$ using (\ref{Mm}) with $\alpha\int_{p-\cdot}^{p_L-v(\cdot)} \rA(s)ds$ replaced by $\alpha \int_{p-\cdot}^p \rA(s)ds-\frac\alpha 6v(\cdot)$. Using  (\ref{v'*}), (\ref{dxiK}), (\ref{dY-m}), and the modified (\ref{dAI}), (\ref{dXK}), (\ref{dA1}), and (\ref{dYK-m}), we find that (\ref{dM/M}) holds here with $\Lambda_0(p_L-v(t),\cdot)$ replaced by $(\frac\kappa 2-3)\tanh_2$. We may write   $M_m=N_m\exp(\cc U)$, where
$$ N_m={A_1^\alpha A_{I,m}^\alpha} \frac{Y_{L,m}}{ Y_m} \exp\Big( (\alpha+\frac {\cc}2)\Big(\int_{p-\cdot}^p\Big( \rA(s)+\frac 1s\Big)ds-\frac{v}6\Big)-\alpha\int_{p-\cdot}^p\frac 1{s}ds\Big);$$ 
$$U=-\frac 16\int_0^\cdot A_S(s)ds+\frac 1{12}v-\frac 12\int_{p-\cdot}^p \Big(\rA(s)+\frac 1{s}\Big)ds.$$


To get estimations on $N_m(t)$, we do some rescaling. Let $\ha p$, $\ha T$, $\check t$, $\ha x_0$, $\ha y_m$, and $\ha\beta$ be as defined in the first paragraph of Section \ref{Section-Rescaling}. Then (\ref{T=p}) holds here.
 From (\ref{rA-ha}), we see that (\ref{int-rA-ha=}) holds if $\ha p_L(t)$ is replaced by $\ha p$ and $p_L-v(\check t)$ is replaced by $p$. 
 Define $\ha\xi(t)$, $\ha q_m(t)$, and $\ha X_m(t)$ using (\ref{haX}); define  $\ha\xi_L$, $\ha q_{L,m}$, and $\ha X_{L,m}$ using (\ref{haXL} ) with the factors $\frac{\ha p_L(t)}\pi$ removed.
Define $\ha g(t,\cdot)$ and $\ha g_{L,W}(t,\cdot)$ using (\ref{ha-g-L}) with the factor $\frac{\ha p_L(t)}\pi$ removed. Then (\ref{conf-g-ha}) holds here and (\ref{conf-g-L-ha}) holds if ``$\St_\pi\sem ((\ha \beta((0,t])+2\ha p\Z)\cup\ha L); \R_\pi$'' is replaced by ``$\ha D_0\sem \ha \beta((0,t]); \ha J_0$'',
where $\ha D_0:=\frac{\ha p}\pi \til D_0$ and  $\ha J_0:=\frac{\ha p}\pi\til J_0$.  Define $\ha W(t,\cdot)$, $\ha A_1(t)$, and $\ha A_{I,m}(t)$ using (\ref{W-ha}) and (\ref{ha-A-def}). Then (\ref{W(q)-ha-m}) still holds,  (\ref{A-ha=}) holds with $\ha p_L(t)$ replaced by $\pi$, and  (\ref{haW=}) should be replaced by $\ha W(t,\cdot):(\ha D_{0,t};\ha J_{0,t})\conf (\St_\pi;\R_\pi)$, 
where  $\ha D_{0,t}:=\ha g(t,\ha D_0)$ and $\ha J_{0,t}:=\ha g(t,\ha J_0)$. 

Let $\ha v(t)=v(\check t)$. Define $\ha Y_m(t)$  using (\ref{Y=X-m-ha}), but define $\ha Y_{L,m}(t)=\ha\Gamma_\infty(\ha v(t),\ha X_{ L,m}(t))$.
Then (\ref{Y-ha=}) holds with  $\ha p_L(t)$ replaced by $\pi$.
Define $\ha N_m$ on $[0,\ha T)$ such that  $$\ha N_m=\Big(\frac{\ha p}\pi\Big)^\alpha\ha A_1^\alpha \ha A_{I,m}^\alpha  {\ha Y_{L,m}}{ \ha Y_m}^{-1}\exp\Big((\alpha+\frac{\cc}2)\Big(\int^{\ha p+ \cdot}_{\ha p} \ha \rA( s)ds-\frac{\ha v}6\Big)\Big).$$
From  the modified (\ref{int-rA-ha=}),  (\ref{A-ha=}) and (\ref{Y-ha=}),  we find that (\ref{N-N}) holds here. 
From (\ref{cc}), (\ref{RA-rA-ha}),  (\ref{sigma-tau}), (\ref{alpha}), (\ref{Gamma-m-ha}), (\ref{ha-Gamma-infty}), and the modified (\ref{Y=X-m-ha}), we see that
$$ \ha N_m=C_p\ha A_1^\alpha \ha A_{I,m}^\alpha  \,\ha\Gamma_q(\ha p+ \cdot,\ha X_m)^{-1}\exp\Big(-\alpha\int_{\ha X_m}^{\ha X_{{L},m}} \tanh_2(s)ds+(\alpha+\frac{\cc}2)\ha\RA(\ha p+\cdot)\Big),$$ 
where $C_p:=(\frac{\ha p}\pi)^\alpha\exp(-(\alpha+\frac{\cc}2)(\ha\RA(\ha p)+\frac{\ha p}6))$.

Let ${\cal E}_m $, $m\in\Z$, be  as in Section \ref{section-est}. Since $\ha \beta(t)$ stays inside $\ha D_0$ before time $\ha T$, we see that $\{\ha T=\infty\}\cap{\cal E}_m=\emptyset$ for $m\in\Z\sem\{m_0\}$.
Suppose that  $\{\ha T=\infty\}\cap{\cal E}_{m_0}$ occurs. An argument using extremal length shows that
$\dist(\{\ha\xi(t),\ha q_m(t)+\pi i\},(\St_\pi\cup\R_\pi)\sem \ha D_{0,t})\to \infty$ as  $t\to\infty$.
Applying Proposition \ref{ThmK} and Proposition \ref{ThmK'}, we find that Proposition \ref{N-to-1-ha} (i) holds here with $m=m_0$ and $C_{p,L}$ replaced by $C_p$.

Let ${\cal P}_m$ denote the family of pairs of disjoint   polygonal crosscuts $(\rho_1,\rho_2)$ in $\ha D_0$ such that, i) for $j=1,2$, the two end points of $\rho_j$ lie on $\R$ and $\R_\pi$, respectively; ii) for $j=1,2$, the line segments of $\rho_j$ are parallel to $x$ or $y$ axes, and all vertices other than the end points have rational coordinates; and iii) $\dist(\rho_1\cup\rho_2,\pa \ha D_0)>0$ and $\rho_1\cup\rho_2$ disconnect $\ha x_0$ and $\ha y_{m}+\pi i$ from $\pa\ha D_0$ in $\St_\pi$. For each $(\rho_1,\rho_2)\in{\cal P}_m$, define $\ha T_{\rho_1,\rho_2}$ to be the biggest time such that $\ha\beta((0,\ha T_{\rho_1,\rho_2}))\cap(\rho_1\cup\rho_2)=\emptyset$. Applying Lemma \ref{comp}, Lemma \ref{comp2} and Lemma \ref{comp3} we find that Proposition \ref{N-to-1-ha} (ii) holds here. We define $\til{\cal P}_m$ as in Section \ref{section-est}. Then Proposition \ref{N-to-1} holds here with with $m=m_0$ and $C_{p,L}$ replaced by $C_p$.

Following the argument of Lemma \ref{BL0} and Lemma \ref{BL}, we can show that $U(t)$ equals the Brownian loop measure of the loops in $\A_p$ that intersect both $L$ and $\beta((0,t))$. Here we use the fact that $-\frac 1{2\pi}\Imm \coth_2(\cdot-x)$ is the Poisson kernel in $\St_\pi$ with the pole at $x\in\R$.

Let $\mu_{L,m}$ denote the distribution of $(\xi(t))$ if $\xi(t)$, $0\le t<T$, is the maximal solution of (\ref{xi-K-m}) with $(\frac\kappa 2-3)\tanh_2$ in place of $\Lambda_0(p_L-v(t),\cdot)$, and satisfies $\xi(0)=x_0$.
Suppose $(\xi(t))$ has distribution $\mu_{L,m_0}$. From (\ref{dxiK}) we conclude that $\beta_L(t)=W_L(\beta(t))$, $0\le t<T$, is a time-change of strip SLE$(\kappa;\kappa-6)$ trace in $\St_\pi$ started from $W_L(e^{ix_0})$ with marked point $W_L(e^{-p+iy_{m_0}})$. Thus, under this distribution, $\beta$ is a time-change of a chordal SLE$(\kappa)$ trace in $\A_p\sem L$ from $z_0=e^{ix_0}$ to $w_0=e^{-p+iy_0}$. Let $\mu_m$ denote the distribution of the maximal solution of $(\ref{xi-p-m})$, or equivalently (\ref{xi-p-m-X}). Using the argument in Section \ref{Section-RN},  Girsanov's theorem, and the modified (\ref{dM/M}) and Proposition \ref{N-to-1} we conclude that for some constant $Z_{m_0}>0$,
\BGE \frac{d\mu_{L,m_0}}{d\mu_{m_0}}=\frac{ {\bf 1}_{\beta\cap L=\emptyset}}{Z_{m_0}} \exp(\cc\mu_{\mbox{loop}}[{\cal L}_{L,p}]).\label{dmuL/dmum0}\EDE

Let $s\in\R$. If  the distribution of $(\xi(t))$ is the $\mu_{\langle s\rangle}$ in Section \ref{section-limit}, then $\beta$ is an annulus SLE$(\kappa;\Lambda_{\langle s\rangle})$ trace in $\A_p$ started from $z_0=e^{ix_0}$ with marked point $w_0=e^{-p+iy_0}$. Since $\{\beta\cap L=\emptyset\}\cap{\cal E}_m=\emptyset$ for $m\in\Z\sem\{m_0\}$, from (\ref{dmu-m/dmu}) we see that (\ref{dmuL/dmum0}) holds with $\mu_{m_0}$ replaced by $\mu_{\langle s\rangle}$ and $Z_{m_0}$  replaced by some other $Z_{\langle s\rangle}>0$. This finishes the sketch of the proof of Theorem \ref{Main-strip}.

\subsection{Multiple SLE crossing an annulus}
Fix $\kappa\in(0,4]$ and $p>0$. Let $n\in\N$ and $n\ge 2$. Let $z_1,\dots,z_n$ be $n$ distinct points that lie on $\TT$ in the counterclockwise direction. Let $w_1,\dots,w_n$ be $n$ distinct points that lie on $\TT_p$ in the counterclockwise direction. Let $\vec z=(z_1,\dots,z_n)$ and $\vec w=(w_1,\dots,w_n)$. Let $\cal G$ denote the set of $(\beta_1,\dots,\beta_n)$ such that each $\beta_j$ is a crosscut in $\A_p$ that connects $z_j$ and $w_j$, and the $n$ curves are mutually disjoint.

\begin{Definition}
  A random $n$-tuple $(\beta_1,\dots,\beta_n)$ with values in $\cal G$ is called a multiple SLE$(\kappa)$ in $\A_p$ from $\vec z$ to $\vec w$ if for any $j\in\{1,\dots,n\}$, conditioned on all other $n-1$ curves, $\beta_j$ is a chordal SLE$(\kappa)$ trace from $z_j$ to $w_j$ that grows in $D_j$, which is the subregion in $\A_p$ bounded by $\beta_{j-1}$ and $\beta_{j+1}$ ($\beta_0=\beta_n$ and $\beta_{n+1}=\beta_1$) that has  $z_j$ and $w_j$ as its boundary points. \label{multiple-SLE}
\end{Definition}

\begin{Theorem}
Let $s_1,\dots,s_n\in\R$. For $j=1,\dots,n$, let $\nu_j$ denote the distribution of the annulus SLE$(\kappa;\Lambda_{\kappa;\langle s_j\rangle})$ trace in $\A_p$ started from $z_j$ with marked point $w_j$. Define a joint distribution $\nu^M$ of $(\beta_1,\dots,\beta_n)$ by
\BGE \frac{d\nu^M}{\prod_{j=1}^n \nu_j}=\frac{{\bf 1}_{{\cal E}_{\disj}}}{Z} \exp\Big(\cc \sum_{s=2}^n\mu_{\mbox{loop}}({\cal L}_{\ge s})\Big),\label{dnu/nu}\EDE
where ${\cal E}_{\disj}$ is the event that $\beta_j$, $1\le j\le n$, are mutually disjoint;  ${\cal L}_{\ge s}$ is the set of loops in $\A_p$ that intersect at least $s$ curves among $\beta_j$, $1\le j\le n$; and $Z>0$ is a  constant. Then $\nu^M$ is the distribution of a multiple SLE$(\kappa)$ in $\A_p$ from $\vec z$ to $\vec w$. \label{multiple}
\end{Theorem}
{\bf Proof.} Suppose for $1\le j\le n$, $\beta_j$ is a crosscut in $\A_p$ connecting $z_j$ with $w_j$.
Fix $j\in\{1,\dots,n\}$. Let ${\cal L}^{j,1}_{\ge s}$ (resp.\ ${\cal L}^{j,0}_{\ge s}$) denotes the set of loops in $\A_p$ that intersect at least $s$ curves among $\beta_k$, $k\ne j$, and intersect (resp.\ do not intersect) $\beta_j$. Then ${\cal L}_{\ge s}={\cal L}^{j,0}_{\ge s}\cup{\cal L}^{j,1}_{\ge s-1}$. Let ${\cal L}^j_{\ge s}={\cal L}^{j,0}_{\ge s}\cup{\cal L}^{j,1}_{\ge s}$. Then ${\cal L}^j_{\ge s}$ depends only on $\beta_k$, $k\ne j$. Since ${\cal L}^{j,0}_{\ge n}=\emptyset$, we have
$$\sum_{s=2}^n\mu_{\mbox{loop}}({\cal L}_{\ge s})= \sum_{s=2}^n\mu_{\mbox{loop}}({\cal L}^{j,0}_{\ge s}) +\sum_{s=1}^{n-1} \mu_{\mbox{loop}}({\cal L}^{j,1}_{\ge s}) = \mu_{\mbox{loop}}({\cal L}^{j,1}_{\ge 1})+ \sum_{s=2}^{n-1}\mu_{\mbox{loop}}({\cal L}^{j}_{\ge s}).$$

Let ${\cal E}^j_{\disj}$ denote the event that $\beta_k$, $k\ne j$, are mutually disjoint. When ${\cal E}^j_{\disj}$ occurs, let  $D_j$ be the simply connected subdomain of $\A_p$ as in Definition \ref{multiple-SLE}. Let $L_j=\A_p\sem D_j$. Then ${\cal E}_{\disj}={\cal E}_{\disj}^j\cap\{\beta_j\cap L_j=\emptyset\}$. Thus, we may rewrite the righthand side of (\ref{dnu/nu}) as
$$\frac{{\bf 1}_{{\cal E}^j_{\disj}}{\bf 1}_{\{\beta_j\cap L_j=\emptyset\}}}{Z}\exp\Big(\cc \sum_{s=2}^{n-1}\mu_{\mbox{loop}}({\cal L}^j_{\ge s})+\cc \mu_{\mbox{loop}}({\cal L}^{j,1}_{\ge 1})\Big)=C_*{\bf 1}_{\{\beta_j\cap L_j=\emptyset\}}\exp(\cc  \mu_{\mbox{loop}}({\cal L}^{j,1}_{\ge 1})),$$
where $C_*=\frac 1Z{{\bf 1}_{{\cal E}^j_{\disj}}}\exp(\cc \sum_{s=2}^{n-1}\mu_{\mbox{loop}}({\cal L}^j_{\ge s}))$ is measurable w.r.t.\ the $\sigma$-algebra generated by $\beta_k$, $k\ne j$. Let $\nu^M_j$ denote the conditional distribution of $\beta_j$ when $(\beta_1,\dots,\beta_n)\sim \nu^M$ and all $\beta_k$ other than $\beta_j$ are given. The above argument shows that the conditional Randon-Nikodym derivative between $\nu^M_j$ and $\nu_j$ is $C_* {\bf 1}_{\{\beta_j\cap L_j=\emptyset\}}\exp(\cc  \mu_{\mbox{loop}}({\cal L}^{j,1}_{\ge 1}))$. Note that ${\cal L}^{j,1}_{\ge 1}$ is the set of all loops in $\A_p$ that intersect both $\beta_j$ and $L_j$. From Theorem \ref{Main-strip} we conclude that $\nu^M_j$ is the distribution of a time-change of a chordal SLE$(\kappa)$ trace in $\A_p\sem L_j=D_j$ from $z_j$ to $w_j$. $\Box$
\vskip 4mm

Choose $x_j,y_j\in\R$ such that $z_j=e^{ix_j}$, $w_j=e^{iy_j-p}$, $1\le j\le n$, $z_1<z_2<\cdots<z_n<z_1+2\pi$, and $w_1<w_2<\cdots<w_n<w_1+2\pi$. For each $m\in\Z$, let ${\cal G}_m$ denote the set of $(\beta_1,\dots,\beta_n)\in\cal G$ such that for each $j$,  $(e^i)^{-1}(\beta_j)$ has a component that connects $x_j$ with $y_j+2m\pi +p i$. Then $\cal G$ is the disjoint union of ${\cal G}_m$'s. Let $\nu^M$ be given by Theorem \ref{multiple}, and let $\nu^M_m=\nu^M[\cdot|{\cal G}_m]$, $m\in\Z$. Then each $\nu^M_m$ is also the distribution of a multiple SLE$(\kappa)$ in $\A_p$ from $\vec z$ to $\vec w$, and the same is true for any convex combination of $\nu^M_m$'s. In fact, the converse is also true.

\begin{Proposition}
  If $\nu$ is the distribution of a multiple SLE$(\kappa)$ in $\A_p$ from $\vec z$ to $\vec w$, then $\nu$ is some convex combination of $\nu^M_m$, $m\in\Z$. \label{unique}
\end{Proposition}
{\bf Proof.}  Define another probability measure $\nu^*$ by
$\frac{d\nu^*}{d\nu}= \frac 1Z \exp\Big(-\cc \sum_{s=2}^n\mu_{\mbox{loop}}({\cal L}_{\ge s})\Big)$,
where $Z>0$ is a normalization constant. From the proof of Theorem \ref{multiple}, we see that, if $(\beta_1,\dots,\beta_n)\sim \nu^*$, then for any $j$, conditioning on  the other $n-1$ curves, $\beta_j$ has the distribution of an annulus SLE$(\kappa;\Lambda_{\langle s_j\rangle})$ trace in $\A_p$ from $z_j$ to $w_j$ conditioned to avoid other curves.

Let ${\cal A}$ denote the set of $( \Omega_1,\dots, \Omega_n)$ such that each $\Omega_j$ is a subdomain of $\A_p$ bounded by two crosscuts crossing $\A_p$, and the $\Omega_j$'s are mutually disjoint. Let ${\cal S}_{\Omega_j}$ denote the event that the curve stays within $\Omega_j$. Let $\mu=\nu^*[\cdot|\prod_{j=1}^n {\cal S}_{\Omega_j}]$. From the property of $\nu^*$, we see that, if $(\beta_1,\dots,\beta_n)\sim \mu$, then for any $j$, conditioning on the other $n-1$ curves, $\beta_j$ has the distribution of an annulus SLE$(\kappa;\Lambda_{\langle s_j\rangle})$ trace in $\A_p$ from $z_j$ to $w_j$ conditioned to stay inside $\Omega_j$. Thus, $\mu=\prod_{j=1}^n \nu_j[\cdot|{\cal S}_{\Omega_j}]$. This implies that $\nu^*=C (\Omega_1,\dots,\Omega_n)\prod_{j=1}^n \nu_j$ on $\prod_{j=1}^n {\cal S}_{\Omega_j}$ for some positive constant $C(\Omega_1,\dots,\Omega_n)$.

Decompose ${\cal A}$ into ${\cal A}_m$, $m\in\Z$, such that ${\cal A}_m$ is the set of all $( \Omega_1,\dots, \Omega_n)\in{\cal A}$ such that there exists $(\beta_1,\dots,\beta_n)\in{\cal G}_m$ with $\beta_j\in\Omega_j$, $1\le j\le n$. Fix $m\in\Z$ and $(\Omega_1,\dots,\Omega_n),(\Omega_1',\dots,\Omega_n')\in{\cal A}_m$. Then $\nu_j({\cal S}_{\Omega_j}\cap{\cal S}_{\Omega_j'})>0$ for each $j$. Thus, $\prod_{j=1}^n {\cal S}_{\Omega_j}\cap  \prod_{j=1}^n {\cal S}_{\Omega_j'}$ is a positive event under $\prod \nu_j$. So we must have $C(\Omega_1,\dots,\Omega_n)=C(\Omega_1',\dots,\Omega_n')$. This means that the function $C( \Omega_1,\dots, \Omega_n)$ is constant, say $C_m$, on each ${\cal A}_m$.
For $m\in\Z$, we may find countably many $(\Omega_1,\dots,\Omega_n)\in{\cal A}_m$ such that the events $\prod_{j=1}^n {\cal S}_{\Omega_j}$ cover ${\cal G}_m$. Thus, $\nu^*=C_m\prod_{j=1}^n \nu_j$ on ${\cal G}_m$ for each $m\in\Z$, which implies that $\nu[\cdot|{\cal G}_m]=\nu^M_m$ for each $m\in\Z$. Since $\nu$ is supported by ${\cal G}=\bigcup_{m\in\Z} {\cal G}_m$, the proof is finished. $\Box$

\vskip 4mm
\no{\bf Remarks.} 
\begin{enumerate}
  \item Theorem \ref{multiple} extends the main result in \cite{LawKoz} which states that, if $\A_p$ is replaced by a simply connected domain $D$, if $z_1,\dots,z_n,w_n,\dots,w_1$ are $2n$ distinct points that lie on $\pa D$ in the counterclockwise direction, if $\nu_j$ is the distribution of a chordal SLE$(\kappa)$ trace in $D$ from $z_j$ to $w_j$, and if $(\beta_1,\dots,\beta_n)$ has joint distribution $\nu^M$ which is defined by (\ref{dnu/nu}), then for any $1\le j\le n$, conditioning on the other $n-1$ curves, $\beta_j$ is a time-change of a chordal SLE$(\kappa)$ trace from $z_j$ to $w_j$ that grows in the component of $D\sem \bigcup_{k\ne j}\beta_k$ whose boundary contains $z_j$ and $w_j$. In fact, for the $(\beta_1,\dots,\beta_n)$ in Theorem \ref{multiple}, if we condition on one of the curves, say $\beta_n$, then the conditional joint distribution of the rest of the curves $\beta_1,\dots,\beta_{n-1}$ agrees with the joint distribution given by \cite{LawKoz} with $D=\A_p\sem\beta_1$.
  \item Since ${\cal G}_m$'s are mutually disjoint, Proposition \ref{unique} implies that for each $m\in\Z$, $\nu^M_m$ does not depend on the choice of $s_1,\dots,s_n$. In fact, if we define another multiple SLE$(\kappa)$ distribution $\nu^{M'}$ using $s_1',\dots,s_n'\in\R$, then there is a constant $Z>0$ such that for each $m\in\Z$, $\frac{d\nu^{M'}}{d\nu^M}=e^{\frac{2\pi m}\kappa(\sum s_j'-\sum s_j')}$ on ${\cal G}_m$. Moreover, since each $\mu_j$ satisfies reversibility, we see that $\nu^M$ and $\nu^M_m$ should also satisfy reversibility.
  \item In the case $n=2$, if we let the inner circle shrink to $0$, it is expected that the two curves tend to the two arms of a two-sided radial SLE$(\kappa)$. The two-sided radial SLE$(\kappa)$ ($\kappa\le 4$) generates two simple curves in $\D$, which connect $0$ with two different points on $\TT$, and intersect only at $0$. The union of the two arms can be understood as a chordal SLE$(\kappa)$ trace connecting the two boundary points, conditioned to pass through $0$. Thus, the knowledge on multiple SLE$(\kappa)$ with $n=2$ can be used to study the microscopic behavior of an SLE trace near a typical point on the trace. 
\end{enumerate}


\begin{thebibliography}{99}
\bibitem{elliptic} K.\ Chandrasekharan. {\it Elliptic functions}. Springer-Verlag Berlin Heidelberg, 1985.
\bibitem{Ahl} Lars V.\ Ahlfors. {\it Conformal invariants: topics
in geometric function theory}. McGraw-Hill Book Co., New York, 1973.
\bibitem{LawKoz} Michael J.\ Kozdron and Gregory F.\ Lawler. The configurational measure on mutually avoiding SLE paths. {\it Universality and Renormalization: From Stochastic Evolution to Renormalization of Quantum Fields}, I. Binder, D.\ Kreimer, ed., Amer. Math. Soc., 199-224.
\bibitem{LawSLE} Gregory F.\ Lawler. {\it Conformally Invariant Processes in the Plane}.
Am.\ Math.\ Soc., Providence, RI, 2005.
\bibitem{LSW-8/3} Gregory F.\ Lawler, Oded Schramm and Wendelin Werner.
Conformal restriction: the chordal case, {\it J.\ Amer.\ Math.\
Soc.}, 16(4): 917-955, 2003.
\bibitem{loop} Gregory F.\ Lawler and Wendelin Werner. The Brownian loop soup. {\it Probab.\ Theory Related Fields},
128(4):565-588, 2004.
\bibitem{RY} Daniel Revuz and Marc Yor. {\it Continuous Martingales
and Brownian Motion}. Springer, Berlin, 1991.
\bibitem{RS-basic} Steffen Rohde and Oded Schramm. Basic properties of
SLE. {\it Ann.\ of Math.}, 161(2):883-924, 2005.
\bibitem{SW} Oded Schramm and David B.\ Wilson. SLE coordinate changes. {\it New York Journal of Mathematics}, 11:659--669, 2005.
\bibitem{thesis} Dapeng Zhan. Random Loewner Chains in Riemann Surfaces. PhD thesis, Caltech, 2004.
\bibitem{Zhan} Dapeng Zhan. Stochastic Loewner evolution in doubly connected domains.
{\it Probab.\ Theory Related Fields}, 129(3):340-380, 2004.
 \bibitem{ann-prop} Dapeng Zhan. Some properties of annulus SLE. {\it Electron.\  J.\  Probab.}, 11, Paper 41:1069-1093, 2006.
\bibitem{LERW} Dapeng Zhan. The Scaling Limits of Planar LERW in Finitely Connected Domains. {\it Ann. Probab.}, 36(2):467-529, 2008.
\bibitem{whole} Dapeng Zhan. Reversibility of whole-plane SLE, arXiv:1004.1865.

\end{thebibliography}
\end{document}